\documentclass[11pt,reqno]{amsart}
\usepackage{amsthm,amsmath}
\usepackage{enumerate}
\usepackage{fancyhdr}
\usepackage{amssymb}
\usepackage{mathrsfs}
\usepackage{verbatim}
\usepackage[body={36cc,51cc}, hmargin={2.8cm,2.8cm}]{geometry}
\usepackage{appendix}
\usepackage{color}
\usepackage{graphicx}
\usepackage{pdfsync}
 \newtheorem{theorem}{Theorem}[section]
 
 \newtheorem{lemma}[theorem]{Lemma}
 \newtheorem{proposition}[theorem]{Proposition}
 \numberwithin{equation}{section}
\allowdisplaybreaks
 \theoremstyle{definition}
 
 \theoremstyle{definition}

 \theoremstyle{remark}

 \newcommand{\bee}{\begin{equation}}
 \newcommand{\ene}{\end{equation}}

\newcommand{\ga}{\gamma}
\newcommand{\les}{\lesssim}

\newcommand{\lt}{\left}
\newcommand{\rt}{\right}

\title[Nonlinear Asymptotic Stability]{Nonlinear asymptotic stability of gravitational hydrostatic equilibrium  for viscous white dwarfs with symmetric perturbations}
\author[Tao Luo, Yan-Lin Wang \& Huihui Zeng]{Tao Luo\textsuperscript{1}, Yan-Lin Wang\textsuperscript{2} \&Huihui Zeng\textsuperscript{3} }
\address[1]{Department of Mathematics, City University of Hong Kong, Kowloon Tong,Hong Kong SAR, China}
\email{taoluo@cityu.edu.hk}
\address[2]{Yau Mathematical Sciences Center, Tsinghua University, Beijing 100084, China}
\email{yanlwang@mail.tsinghua.edu.cn}
\address[3]{Department of Mathematical Sciences, Tsinghua University, Beijing 100084, China}
\email{hhzeng@mail.tsinghua.edu.cn}
\begin{document}

\begin{abstract} We prove the nonlinear asymptotic stability of the gravitational hydrostatic equilibrium for the general equation of state of pressure-density relation in the framework of vacuum free boundary problem of spherically symmetric compressible Navier-Stokes-Poisson equations in three dimensions.
The results apply to white dwarfs and polytropes for all $\gamma>4/3$ including the case of $\gamma\geq 2$ which was not addressed in previous literature.
Detailed decay rates of perturbations are given.

\vspace{0.5cm}

\end{abstract}

\maketitle

\section{Introduction}
\subsection{Motivation}
In the fundamental hydrodynamical setting (cf. \cite{chandr}), the gravitational hydrostatic equilibrium of a star satisfies the following  equation:
 \begin{equation}\label{1}p(\bar \rho)_r+ 4\pi\mathcal{G}\bar \rho r^{-2}\int_0^r  s^2 \bar\rho(s) ds=0 \ \ {\rm for} \ \ r>0, \end{equation}
where $\bar\rho(r)$ is the density, $p(\bar \rho)$ is the pressure and $\mathcal{G}$ is the gravitational constant. We assume that the pressure satisfies
 \begin{subequations}\label{pcon1} \begin{align}
&p\in C^1\left([0, \infty)\right)\cap   C^2\left((0, \infty)\right), \ \  p'(s)>0  \ {\rm for} \ s>0,   \label{pcon1-a}\\
& \lim_{s\to 0+}s^{-4/3}p(s)=0,  \ \ \lim_{s\to \infty} s^{-4/3}p(s)=\mathcal{K} \label{pcon1-b}
\end{align}\end{subequations}
for either $ \mathcal{K}\in (0, \infty)$ or  $\mathcal{K}=\infty$. There are two typical stars, {\em polytropes} satisfying
\begin{equation}\label{polytropic}
p(\rho)=\kappa \rho^\gamma  \ \ \textrm{with positive constnats $\kappa$ and $\gamma$}, \end{equation}
and {\em white dwarfs} satisfying \eqref{pcon1} with $ \mathcal{K}\in (0, \infty)$.  For a  white dwarf, the gravity of which is balanced by electron degeneracy pressure (not dependent on the temperature of the white dwarf), and the pressure $p$ as a function of density $\rho$ is given through the following implicit relation (cf. \cite{chandr}):
$$
p(x)=\Gamma_1 (x(2x^2-3)(x^2+1)^{1/2}+3\sinh^{-1} x) \ \ {\rm and} \ \ \rho(x)=\Gamma_2 x^3
$$
for $x\ge 0$,
and satisfies the following asymptotics (cf. \cite{chandr}):
\begin{equation*}p(\rho)=\lt\{\begin{split}
&c_1\rho^{4/3}-c_2\rho^{2/3}+\cdots, \quad \rho\rightarrow \infty,\\
&d_1\rho^{5/3}-d_2\rho^{7/3}+O( \rho^3), \quad \rho\rightarrow 0,\\
\end{split}\rt.
\end{equation*}
where $\Gamma_1$, $\Gamma_2$, $c_1, c_2, d_1, d_2$ are positive physical constants.

The existence, uniqueness and properties
of gravitational hydrostatic equilibria with the prescribed total mass are summarized (cf. \cite{ab, liebyau, LSmoller, Rein}) as follows:
\begin{proposition}\label{prop1}
Suppose that the pressure satisfies \eqref{pcon1}.
There exists a constant $M_c$ satisfying $M_c\in (0, \infty)$ if $\mathcal{K}\in (0, \infty)$ and $M_c=\infty$ if  $\mathcal{K}=\infty$
such that if $M\in (0,M_c)$, then there exist a unique $\bar R>0$ and a unique solution $\bar \rho(r)$ to problem \eqref{1} satisfying that $\bar\rho(r)=0$ for $r\ge \bar R$, $\bar\rho_r(0)=0$,  $\bar\rho_r(r)\in (-\infty, 0)$ for $r\in(0,\bar R)$ and
 $4\pi \int_0^{\bar R} s^2 \bar\rho(s)ds=M$.
In this case,  $\bar \rho$ is a minimizer of the energy functional
\begin{equation}\label{energy}
E(\rho)=\int_{\mathbb{R}^3} \mathcal{A}(\rho) dx + 2^{-1}\int_{\mathbb{R}^3}\rho \Psi  dx \end{equation}
satisfying the total mass constraint
$\int_{\mathbb{R}^3} \rho dx=M$ in the class of nonnegative functions,
where
$$\mathcal{A} (\rho)=\rho\int_0^\rho s^{-2}p(s) ds \ \ {\rm and} \ \  \Psi(\rho)(x)=-\mathcal{G}\int_{\mathbb{R}^3} \frac{\rho (y)}{|x-y|}dy. $$
The first term on the right-hand-side of \eqref{energy} represents the potential energy, and the second stands for the  gravitational energy.
\end{proposition}

The stellar material  for a white dwarf by virtue of Pauli's exclusion principle and the uncertainty principle  exercises a ground state pressure, which is the sole local balance for the gravitational force in a non-rotating white dwarf,  since no more nuclear fuel to burn to supply additional thermal and radiation pressure gradients.
In the famous Chandraskhar theory, gravitational collapse occurs when the total mass exceeds a critical mass, the ``Chandraskhar Limit"  (cf. \cite{chandr, liebyau, wein}), which is $M_c$ in Proposition \ref{prop1}.
The problem of the gravitational collapse for white dwarfs  was
formulated in \cite{4,8,16} leading  to the equation for the density, equation  \eqref{1},  called the ``Chandrasekhar equation"  in \cite{liebyau}.
The gravitational collapse was predicted in \cite{chandr,4} by this equation at some critical mass, which was  verified in \cite{liebyau} as the limit
of quantum mechanics.
It was also proved in \cite{liebyau} that, for a white dwarf, there exists $M_c\in(0,\infty)$ such that for each $M\in (0,M_c)$,
there exists a unique radially symmetric decreasing minimizer $\rho_M$ for the energy functional \eqref{energy} with the total mass constraint
$\int_{\mathbb{R}^3} \rho_M dx=M$, which satisfies the Euler-Lagrangian equation for some Lagrangian multiplier $\mu_M$:
$$ \mathcal{A}'(\rho_M(x))=\{-\Psi(\rho_M)(x)-\mu_M\}_+, $$
where $\{f(x)\}_+=\max\{f(x), 0\}$, and $\mu_M\to \infty$ as $M\to M_c$.
When the total mass $M$ is less than the ``Chandrasekhar limit'' $M_c$,  the dynamical stability of the gravitational hydrostatic equilibrium (the stationary solution) becomes a very interesting and important problem for white dwarfs.

In this direction,
a conditional nonlinear Lyapunov type
stability theory was established in \cite{LSmoller} by use of a variational approach,  assuming the existence of global solutions of the Cauchy problem for  the three-dimensional compressible Euler-Poisson equations. (The same type of nonlinear stability results can be found in \cite{Rein} for gaseous stars with $\mathcal{K}=\infty$ in \eqref{pcon1-b},  and for  rotating
stars in \cite{LSmoller, LSmoller2}.)
However, rigorous mathematical justification of the nonlinear large time asymptotic stability of the gravitational hydrostatic equilibrium (the stationary solution)  for white dwarfs is missing, even on the linear level.
This is also true for polytropes with $\gamma\ge 2$, though this type of stability was proved in \cite{LXZ2, LZeng} for viscous polytropic stars of \eqref{polytropic}  with $4/3<\gamma<2$.
The main purpose of this article is to address these issues. Precisely, we consider the vacuum free boundary problem for compressible Navier-Stokes-Poisson equations  with a very general equation of state including white dwarfs and polytropes with all $\gamma>4/3$ in three dimensions with spherical symmetry, and prove that if the initial data are small perturbations of the gravitational hydrostatic equilibrium with the same total mass, then there exists a unique strong solution to the vacuum free boundary problem which converges to the gravitational hydrostatic equilibrium as time goes to infinity, with the detailed convergence rates.

The vacuum boundary considered in the present work  satisfies the {\em physical vacuum} condition  that $\sqrt{p'(\rho)}$  is $C^{1/2}$-H\"{o}lder continuous near vacuum states.
We remark that the physical vacuum has been attracting much attention and interests in the study of the compressible fluids (cf. \cite{CHWW}-\cite{CS2},
\cite{DuanQ}-\cite{FZ2}, \cite{GL}-\cite{hlz}, \cite{jang3}-\cite{jang6}, \cite{tpliudamping}-\cite{LY}, \cite{LXZ1}-\cite{LZ}, \cite{YT}-\cite{Zeng2}).
Besides for white dwarfs, the results obtained in this paper
also apply  to viscous polytropic stars of \eqref{polytropic}   with $\gamma>4/3$, as mentioned earlier.
The value of $4/3$ is critical  for the stability of polytropic stars.
In fact, the energy  (cf. \cite{chandr, wein}) for a polytrope with $\gamma>6/5$ on the gravitational hydrostatic equilibrium is given by
$$E_{\gamma}=-\frac{3\gamma-4}{5\gamma-6}\frac{\mathcal{G} M^2}{\bar R}, $$
where $M$ is the total mass and $\bar R$ is the radius of the star. It should be noted that for a white dwarf, the exponent  $4/3$ is the critical  case for
the pressure $p(\rho)$ as $\rho\to \infty$, which is in sharp contrast to the case considered in \cite{LXZ2, LZeng}.

The present work does not cover the case for supermassive stars, which are polytropes with $\gamma=4/3$  and  supported by the pressure of radiation, also called  extreme relativistic degenerate stars (cf. \cite{Shapiro, wein}).
Indeed, for the non-rotating gravitational hydrostatic equilibrium of supermassive stars, the total energy $E_{4/3}=0$,  quoting (\cite{wein}, p. 327)
``the polytrope with $\gamma= 4/3$ is trembling between stability
and instability", and it is remarked that the General Relativity is needed to settle this stability problem.

\subsection{Formulation of the problem and main theorems}
The evolution of the free boundary which is the interface between the gas and vacuum for a viscous gaseous star with spherical symmetry can be modeled by the following free boundary problem:
\begin{subequations}\label{sphereqn}
\begin{align}
&(r^2 \rho)_t+(r^2\rho u)_r=0  & \text{in}\quad(0, R(t)), \label{sphereqn-1}\\
&\rho(u_t+uu_r)+p_r+ \frac{4\pi \mathcal{G} \rho}{ r^{2}}\int_0^r s^2 \rho(s, t) ds =\nu\lt[\frac{(r^2 u)_r}{r^2}\rt]_r & \text{in}\quad(0, R(t)), \label{sphereqn-2}\\
&\rho >0  &  \text{in} \quad[0, R(t)),\\
&\rho (R(t),t)=0 ,
\ \ u(0,t)=0, & \\
&\frac{4}{3}\nu_1 \lt(u_r-\frac{u}{r}\rt) +\nu_2\lt(u_r+2\frac{u}{r}\rt)=0 &{\rm for} \ r=R(t),  \\
& \dot R(t)=u(R(t), t) \ \  {\rm with} \ \ R(0)=R_0, &  \\
&(\rho, u)=(\rho_0, u_0)  & \text{on}\quad(0, R_0).
\end{align}
\end{subequations}
Here $(r,t)$,  $ \rho$,  $u $ and $p$ denote, respectively, the space and time variable, density, velocity and the pressure; the positive constants $\nu_1$, $\nu_2$  and $\mathcal{G}$ represent, respectively, the shear viscosity, the bulk viscosity and the gravitational constant, and $\nu=4\nu_1/3+\nu_2$; and
$R(t)$ is the free boundary of the gaseous star.   This formulation can be found in \cite{LXZ2}.

Let $(\rho,u,R(t))=(\bar\rho(r), 0, \bar R)$ be the stationary solution to \eqref{sphereqn} with the total mass  $M=4\pi\int_0^{\bar R} s^2 \bar\rho(s)ds\in (0,M_c)$  as stated in Proposition \ref{prop1}.  Then it holds that for $r\in  (0, \bar R)$,
\begin{subequations}\label{yo.1}
\begin{align}
&(p(\bar\rho))_r=-r \bar\rho \phi, \ \   \bar\rho_r\in (-\infty,0), \label{pxphi}\\
& \bar\rho_0:= \bar\rho(0) > \bar\rho(r)> \bar\rho(\bar R)= 0,
\label{ibm}
\end{align}\end{subequations}
where
$
\phi(r)= {4\pi \mathcal{G}}{r^{-3}}\int_0^r  s^2 \bar\rho(s)ds \in \left[ {M\mathcal{G}}{\bar R^{-3}}, \ {4\pi\mathcal{G} \bar\rho_0}/{3}\right]
$.
Clearly, it holds that for $r\in (0, \bar R)$,
\begin{align}\label{irhor}
2^{-1} {\bar R}^{-2} M \mathcal{G} (\bar R-r)\le i(\bar\rho(r))\le  ({4\pi \mathcal{G}\bar\rho_0\bar R}/{3}) (\bar R-r),
\end{align}
where the enthalpy $i$ is defined by
$i(s)=\int_0^s \tau^{-1}{p'(\tau)}d\tau$ for $s\ge 0$.
In order to capture the behavior \eqref{irhor} of the stationary solutions near the vacuum boundary, we assume that the initial density satisfies the following condition:
\begin{equation}\label{inidens}
\begin{split}
&\rho_0(r)>0  \ \ {\rm for} \ \  r\in [0,R_0),   \ \  M=4\pi\int_0^{R_0} s^2 \rho_0(s)ds \in (0, \infty), \\
& \rho_0(R_0)=0, \ \  \lt(i(\rho_0)\rt)_r \in (-\infty, 0) \ \  {\rm at}\ \  r=R_0.
\end{split}\end{equation}
That is,  $i(\rho_0(r))$ is equivalent to $R_0-r$ as $r$ is close to $R_0$.
Under the additional assumption:
\begin{align}
\frac{4}{3}\le \inf_{0<s\le \bar\rho_0}\frac{sp'(s)}{p(s)}
\ \ {\rm and} \ \
\bar{\gamma}:=\lim_{s\to 0+}\frac{sp'(s)}{p(s)} \in \lt(\frac{4}{3},\ \infty\rt),
\label{pcon2}
\end{align}
we  prove in this article that the global existence of the strong solution to the vacuum free boundary problem \eqref{sphereqn} and \eqref{inidens} with a general class of pressure laws including white dwarfs and polytropes with $\gamma>4/3$, and the nonlinear asymptotic stability of the gravitational hydrostatic equilibria, extending the results obtained in \cite{LXZ2,LZeng} for polytropes with $4/3<\gamma<2$.

We use the Lagrangian coordinates to fix the domain, as did in \cite{CS1, CS2}.  Let $x\in I=(0, \bar R)$ be the reference variable,
we define the Lagrangian variable $r(x, t)$
by
\begin{align*}
&r_t(x, t)=u(r(x, t), t) \ \text{for}\ t>0,\ \
r(x, 0)=r_0(x) \  {\rm for} \ x\in I,
\end{align*}
where $r_0$ is a diffeomorphism from $\bar I$ to $[0, R_0]$ defined by $\int_0^{r_0(x)}s^2 \rho_0(s)ds
=\int_0^x s^2 \bar\rho(s) ds$. Clearly,
$r_0^2(x)\rho_0(r_0(x))r_0'(x)=x^2 \bar\rho(x)$ for $x\in \bar I$.
Then, it follows from \eqref{sphereqn-1} that
\bee\label{r0x}
\rho(r(x, t), t)=\frac{ r_0^2(x)\rho_0(r_0(x))r_{0x}(x)}{r^2(x, t)r_x(x, t)} =\frac{x^2\bar\rho(x)}{r^2(x, t)r_x(x, t)}.
\end{equation}
We set
$ v(x, t)=u(r(x, t), t)
$ and use \eqref{r0x} to rewrite  problem \eqref{sphereqn} as  the following initial boundary value problem:
\begin{subequations}\label{fixedp}
\begin{align}
&\bar\rho\frac{x^2}{r^2}v_t
+\lt[p\lt(\frac{x^2\bar\rho}{r^2r_x}\rt)
\rt]_x+4\pi \mathcal{G}\frac{x^2}{r^4}\bar\rho\int_0^x y^2\bar\rho dy=\nu \lt(\frac{(r^2 v)_x}{r^2r_x}\rt)_x & \text{in} &\ I\times(0,T], \label{fixedp-a}\\
&v(0, t)=0 \ {\rm and} \ \mathscr{B}(\bar R, t)=0 & \text{on} & \  (0, T],\label{fixedp-b}\\
&(r, v)(x, 0)=(r_0(x), u_0(r_0(x))) & \text{on}& \ I,\label{fixedp-c}
\end{align}\end{subequations}\
where
\begin{equation}\label{mb}
\mathscr{B}:=\frac{4}{3}\nu_1\lt(\frac{v_x}{r_x}-\frac{v}{r}\rt)
+\nu_2\lt(\frac{v_x}{r_x}+2\frac{v}{r}\rt)
=\frac{4}{3}\nu_1\frac{r}{r_x}\lt(\frac{v}{r}\rt)_x
+\nu_2\frac{(r^2v)_x}{r_xr^2}.
\end{equation}
In the setting, the moving vacuum boundary for problem \eqref{sphereqn} is the given by
$
R(t)=r(\bar R, t)
$.

The main result of this paper is stated as follows. We set
\begin{equation}\label{E}\begin{split}
\mathfrak{E}(t)=& \lt\|(r_x-1, v_x)(\cdot, t)\rt\|^2_{L^\infty(I)}
+\lt\|\bar\rho^{1/2}v_t(\cdot, t)\rt\|_{L^2(I)}^2 \\
& +\lt\|\bar\rho^{-1/2} p(\bar \rho)\lt((r/x)_x,\ r_{xx}\rt)(\cdot, t)\rt\|_{L^2(I)}^2.
\end{split}\end{equation}

\begin{theorem} \label{theo1}
Assume that $p$ satisfies  \eqref{pcon1} and \eqref{pcon2}, $\bar \rho$ is a stationary solution  with the total mass $M=4\pi \int_0^{\bar R} s^2 \bar\rho(s)ds\in (0, M_c)$ as stated in Proposition \ref{prop1},
the initial density $\rho_0$ satisfies \eqref{inidens} and $\int_0^{R_0} s^2 \rho_0(s)ds=\int_0^{\bar R} s^2 \bar\rho(s)ds$,
and the compatibility conditions $v_0(0)=0$ and $\mathscr{B}(\bar R, 0)=0$
 hold.
There exists a constant $\bar\delta>0$ such that if  $\mathfrak{E}(0)\le \bar\delta$, then problem \eqref{fixedp} admits a   solution in $I\times[0,\infty)$ with
$$
\mathfrak{E}(t)\leq C \mathfrak{E}(0),\ \ t\geq 0,
$$
for some positive constant $C$ independent of $t$.
\end{theorem}

As a corollary of Theorem \ref{theo1}, we have the nonlinear asymptotic stability of the stationary solution to problem \eqref{sphereqn} as follows:

\begin{theorem}\label{theo2}
Suppose that the assumptions in Theorem \ref{theo1} holds.
There exists a constant $\bar\delta>0$ such that if  $\mathfrak{E}(0)\le \bar\delta$, then problem \eqref{sphereqn} admits a  global solution $(\rho,u, R(t))$ for $t\in [0,\infty)$ satisfying $R\in W^{1,\infty}([0,\infty))$ and the following estimates: for any $\theta\in \lt(0, 1-5/(4\bar\gamma)\rt]$ and $l\in (0, \bar R)$, there exist positive constants $c$,  $C$, $C(\theta^{-1})$ and $C(\theta^{-1},l^{-1})$ independent of $x$ and $t$ such that for $t\in[0, \infty)$,
\begin{subequations}\begin{align}
&c^{-1}\lt(R(t)-r(x,t)\rt) \le  i\lt( \rho(r(x,t), t) \rt) \le c \lt(R(t)-r(x,t)\rt), \ \  x\in I,
\label{thest0}\\
&  \bar\rho^{-2}(x) |\rho(r(x, t), t)- \bar\rho(x)|^2
\le C \mathfrak{E}(0),\ \ x\in I, \label{thest1}\\
& \Upsilon(1+t)^{2- {2}/{\bar\gamma}-3\theta/2}  |r(x,t)-x|^2 + (1+t)^{2\zeta-1/2} \lt\{ (1-\Upsilon) |r(x,t)-x|^2 \rt.\notag\\
&\ \ \lt. +u^2(r(x,t),t)\rt\}
+(1+t)^{2\zeta-1}\lt\{|u_r(r(x, t), t)|^2+|r^{-1}u(r(x, t), t)|^2\rt.\notag\\
&\ \ \lt. +\Upsilon \bar\rho^{(3\bar\gamma-6)/2}(x) |\rho(r(x, t), t)- \bar\rho(x)|^2\rt\}\le C(\theta^{-1}) \mathfrak{E}(0)  ,\ \ x\in I, \label{thest2}\\
&(1+t)^{\min\lt\{\frac{4-\bar\gamma}{2\bar\gamma}
-\theta\lt(\frac{\bar\gamma-1}{2\bar\gamma}
+\frac{4-\bar\gamma-\theta(\bar\gamma-1)}
{4\bar\gamma-2}\rt), \ 2\zeta-1 \rt\}} |\rho(r(x, t), t)- \bar\rho(x)|^2\notag\\
&\ \
\le C(\theta^{-1}) \mathfrak{E}(0), \ \  x\in I, \ \ {\rm if} \  2 < \bar\gamma<4 \ {\rm and}\  \theta<({4-\bar\gamma})/({\bar\gamma-1}),
\label{thest3}\\
& (1+t)^{2\zeta-1}\lt(|r_x(x, t)-1|^2 +|x^{-1}r(x, t)-1|^2\rt)\notag\\
&\ \ \le C(\theta^{-1},l^{-1})\mathfrak{E}(0) , \ \  x\in [0, \bar R -l].\label{thest4}
\end{align}\end{subequations}
Here $\zeta=1- {1}/({2\bar\gamma})
- \theta/2$, and $\Upsilon=1 $ for $\bar\gamma\le 2$ and $\Upsilon=0 $ for $\bar\gamma>  2$.
\end{theorem}

Due to the high degeneracy of   system \eqref{sphereqn} caused by the behaviour \eqref{irhor} near the vacuum boundary, as pointed out in \cite{LXZ2, LZeng} for viscous polytropic stars of \eqref{polytropic} with $4/3<\gamma<2$,  it is extremely challenging to establish  the global-in-time regularity of higher-order derivatives of solutions uniformly up to the vacuum boundary.
Compared with the polytropes considered in \cite{LXZ2, LZeng}, it is way intricate  to build up the higher-order uniform-in-time regularity up to the vacuum  boundary for the general pressure law considered in this paper.

\subsection{Review of related results}
Besides the results already mentioned earlier, we review some related results concerning the mathematical study of white dwarfs and
the stability of gaseous stars.
The existence of steady solutions of rotating and non-rotating
white dwarfs was proved in \cite{ab}, which was extended in \cite{lions,FT}. With the constant angular velocity, the existence and properties such as the bounds of diameters and number of connected components of rotating white dwarfs can be found in \cite{CL,lyy}. The related properties of free boundaries for steady solutions are referred to \cite{CF}.

In \cite{linss}, it was investigated that the linearized stability and instability of gaseous stars including polytropes with $\gamma\in (4/3, 2)$ and $\gamma\in (6/5, 4/3)$,  respectively; concerning that
a turning point principle was proved in a recent work \cite{linzeng} modeled by Euler-Poisson equations.
The nonlinear instability was proved for polytropes with $\gamma=6/5$ in \cite{jang1} and $6/5<\gamma<4/3$ in \cite{jang4}, respectively, in the framework of   Euler-Poisson equations; and for polytropes with $6/5<\gamma<4/3$ in \cite{jang6} in the framework of   Navier-Stokes-Poisson equations. Instability result of Euler-Poisson equations for polytropes with $\gamma=4/3$ was identified  in \cite{DLYY}, since a small perturbation can cause part of the mass to go off to infinity.
For  Euler-Possion equations modelling polytropic stars with $p(\rho)=\rho^{\gamma}$,  a class of global expanding solutions was constructed in \cite{HJ2} for $\gamma=4/3$ and in \cite{HJ1} for other values of $\gamma$, and the continued gravitational collapse was demonstrated in \cite{GHJ2} for $1<\gamma<4/3$  , (see also further study in \cite {HY}) , and  the Larson-Penston self-similar gravitational collapse was proved in \cite{GHJ2} for $\gamma=1$.  The global existence of weak spherically symmetric solutions in multi-dimensions was proved in \cite{CHWY} for certain range of $\gamma$.

The rest of this paper is devoted to the proofs of  Theorems \ref{theo1} and \ref{theo2} by establishing the uniform-in-time higher-order regularity up to the vacuum boundary, and the detailed decay rates of perturbations.

\section{Proof of Theorems \ref{theo1} and \ref{theo2} }

The local existence of  solutions  to problem \eqref{fixedp} in the functional space $\mathfrak{E}(t)$ can be obtained easily by combining the approximation techniques in \cite{jang3,LXZ2} and the following a priori estimates stated in Proposition \ref{22.1.4}.
Indeed, the a priori estimates obtained
in Proposition \ref{22.1.4} are  sufficient for the local existence theory, at least for small $\mathfrak{E}(0)$.
Hence, Theorem \ref{theo1} follows from a standard continuation argument with the following estimates:

\begin{proposition}\label{22.1.4}
Suppose that the assumptions in Theorem \ref{theo1} holds.
Let $v$ be a   solution to problem \eqref{fixedp} in the time interval $[0, T]$
with $r(x, t)=r_0(x)+\int_0^t v(x,\tau) d\tau$
satisfying $\mathfrak{E}(t)<\infty$ on $[0, T]$ and the a priori assumption:
\begin{subequations}\label{21apri}
\begin{align}
&\lt\|\lt(r_x-1, \ {r}/{x}-1\rt)
(\cdot, t) \rt\|_{L^\infty}\le \varepsilon_0, \ \   t\in [0,T],\label{priori1}\\
& \lt\|\lt(v_x, \ {v}/{x}\rt)
(\cdot, t) \rt\|_{L^\infty}\le 1 , \ \  t\in [0,T] \label{priori2}
\end{align}
\end{subequations}
for some suitably small fixed number $\varepsilon_0\in (0, 1/8]$ independent of $t$.
Then there exist positive constants
$C$, $C(\theta^{-1})$ and $C(\theta^{-1},l^{-1})$ independent of $t$ such that for $t\in [0,T]$,
\begin{subequations}\begin{align}
&  \mathfrak{E}(t)+ \lt\|  ({r}/{x}-1, v/x)(\cdot,t)\rt\|_{L^\infty}^2
\le C  \mathfrak{E}(0), \label{22.1.11-a}\\
&(1+t)^{2\zeta} \lt\| (x^{1/2} v )(\cdot,t)\rt\|^2_{L^\infty}
 + (1+t)^{2\zeta-1/2} \lt\{ \lt\|  v (\cdot,t)\rt\|^2_{L^\infty} \rt.\notag \\
 &\lt. + (1-\Upsilon) \lt\|(r-x)(\cdot,t)\rt\|^2_{L^\infty} \rt\}+ \Upsilon(1+t)^{2- {2}/{\bar\gamma}-3\theta/2}  \lt\|(r-x)(\cdot,t)\rt\|^2_{L^\infty} \notag\\
 & +  (1+t)^\zeta \lt\| (x^{3/2} v_x)(\cdot,t)\rt\|^2_{L^\infty}
+(1+t)^{2\zeta-1}\lt\{\lt\|(v_x,v/x) (\cdot, t)\rt\|_{L^\infty}^2  \rt.\notag\\
&  \lt. +\Upsilon \lt\| \lt(\bar\rho^{(3\bar\gamma-2)/4}  \lt(x^2/(r^2 r_x)-1\rt)\rt)(\cdot,t)\rt\|_{L^\infty}^2 \rt\}\le C(\theta^{-1}) \mathfrak{E}(0),  \label{22.1.11-b}\\
&(1+t)^{\min\lt\{\frac{4-\bar\gamma}{2\bar\gamma}
-\theta\lt(\frac{\bar\gamma-1}{2\bar\gamma}
+\frac{4-\bar\gamma-\theta(\bar\gamma-1)}
{4\bar\gamma-2}\rt), \ 2\zeta-1 \rt\}} \lt\| \lt(\bar\rho \lt(x^2/(r^2 r_x)-1\rt)\rt)(\cdot,t)\rt\|_{L^\infty}^2
\notag\\
&
\le C(\theta^{-1}) \mathfrak{E}(0), \ \ {\rm if} \  2 < \bar\gamma<4 \ {\rm and}\  \theta<({4-\bar\gamma})/({\bar\gamma-1})
,\label{22.1.11-c}\\
& (1+t)^{2\zeta-1}\lt(\lt\|\lt(r_{xx}, (r/x)_x, v_{xx}, (v/x)_x\rt)(\cdot, t)\rt\|_{L^2([0,\bar R-l ])}^2\rt. \notag\\
&  \lt.+\lt\|(r_x-1,r/x) (\cdot, t)\rt\|_{L^\infty([0,\bar R-l ])}^2 \rt) \le C(\theta^{-1},l^{-1})\mathfrak{E}(0) .\label{22.1.11-d}
\end{align}\end{subequations}
Here $L^\infty=L^\infty(I)$ in \eqref{21apri} and \eqref{22.1.11-a}-\eqref{22.1.11-c},
$\zeta=1- {1}/({2\bar\gamma})
- \theta/2$, and
\begin{align}
\Upsilon=1  \ \ {\rm for} \  \  \bar\gamma\le 2,
 \ \ \Upsilon=0 \ \  {\rm for} \  \  \bar\gamma > 2. \label{22.1.11-1}
\end{align}
Moreover, there exists a positive constant $c$ independent of $x$ and $t$ such that for
$(x,t)\in I\times [0,T]$,
\begin{align}c^{-1}(\bar R -x) \le i\lt( \lt({\bar\rho x^2 }/({r^2 r_x }) \rt)(x,t)\rt)\le c(\bar R -x). \label{22.1.11-e}
\end{align}
\end{proposition}
Indeed, \eqref{22.1.11-a} follows from \eqref{22.1.5}; \eqref{22.1.11-b} from
\eqref{12.21-da} and
\eqref{12.28-9};  \eqref{22.1.11-c} from
\eqref{20220111}; \eqref{22.1.11-d} from
\eqref{22.1.6-1}; and \eqref{22.1.11-e} from \eqref{har-1} and \eqref{til-3}.

To prove Theorem \ref{theo2}, we recall that
\begin{align*}
\rho(r(x, t), t)= \frac{\bar\rho(x) x^2}{r^2(x, t)r_x(x, t)},\
 u (r(x, t), t)=v(x,t) \  {\rm and} \  R(t)=r(\bar R, t),
 \end{align*}
and notice that
$$
(1/2) x \le r(x,t) \le (3/2) x \ {\rm and}  \  1/2\le r_x(x,t) \le 3/2,
$$
which is due to  \eqref{22.1.11-a} and the smallness of $\mathfrak{E}(0)$. \eqref{thest0} follows from \eqref{22.1.11-e}  and
$$\frac{1}{2}(\bar R -x)\le
R(t)-r(x,t)=\int_x^{\bar R} r_y(y,t)dy\le \frac{3}{2}(\bar R -x);$$
\eqref{thest1} from \eqref{22.1.11-a};
\eqref{thest2} from \eqref{22.1.11-b}, $u/r=(v/x)/(r/x)$ and $u_r=v_x/r_x$;
\eqref{thest3} from \eqref{22.1.11-c}; and
\eqref{thest4} from \eqref{22.1.11-d}.

The rest of this section devotes to proving Proposition \ref{22.1.4}, which consists of the preliminaries in Section \ref{22.1.4-1}, lower-order estimates in Section \ref{22.1.4-2}, and higher-order estimates in Section \ref{22.1.4-3}.

\subsection{Preliminaries}\label{22.1.4-1}
First, it should be noted that
$$
1/2 \le r_x \le 3/2 \ \ {\rm and} \ \ 1/2 \le r/x \le 3/2 \ \ {\rm for}\ \  (x, t)\in I\times [0,T],
$$
which is due to \eqref{priori1}.
We list some properties for functions $p$,  which will be useful in the energy estimates.
It follows from \eqref{pcon1} and \eqref{pcon2}  that
\begin{align}
& p(0)=0, \label{d-0} \\
& \lim_{s\to 0+}\frac{s^{a+1}}{p(s)}\frac{d}{ds} \left(\frac{p(s)}{s^a}\right)
=\lim_{s\to 0+} \frac{sp'(s)}{p(s)}-a=\bar{\gamma}-a \label{d-1}
\end{align}
for any $a>0$.
If $a<\bar{\gamma}$, then there exists a positive  constant $\delta_1$ such that
$\frac{d}{ds} \left(s^{-a} {p(s)} \right)>0$ for
$s\in (0, \delta_1]$ and $\lim_{s\to 0+}  s^{-a} {p(s)} \le {\delta_1^{-a}}{p(\delta_1)}<\infty$.
This, together with the L'Hopital rule, \eqref{d-0} and \eqref{pcon2}, means
\begin{align}
&\lim_{s\to 0+} \frac{p(s)}{s^a}
=\frac{1}{a} \lim_{s\to 0+} \frac{sp'(s)}{p(s)}  \frac{p(s)}{s^a}=\frac{\bar{\gamma}}{a} \lim_{s\to 0+} \frac{p(s)}{s^a}.\label{d-2}
\end{align}
So,
\begin{align}
\lim_{s\to 0+} {s^{-a}} {p(s)}=0  \ \ {\rm for} \ \   a <\bar{\gamma}. \label{d-3}
\end{align}
In particular,
\begin{subequations}\label{d-00}
\begin{align}
&p'(0)=\lim_{s\to 0+} s^{-1}p(s)=0,\\
& \lim_{s\to 0+} s^{-1/6}p'(s)=\bar\gamma \lim_{s\to 0+} s^{-7/6}p(s)=0.
\end{align}\end{subequations}
Similarly, if $a>\bar{\gamma}$, there exists a positive  constant $\delta_2$ such that
$\lim_{s\to 0+} s^{-a} {p(s)}  \ge {\delta_2^{-a}}{p(\delta_2)}>0$.
Suppose that
$\lim_{s\to 0+} s^{-a}{p(s)}<\infty$,
then it yields from \eqref{d-2} that
$\lim_{s\to 0+} s^{-a}{p(s)}=0$.
This is a contradiction. Therefore,
\begin{align}
\lim_{s\to 0+} {s^{-a}} {p(s)}=\infty  \ \ {\rm for} \ \  a >\bar{\gamma}. \label{d-5}
\end{align}
Due to \eqref{pcon1},  \eqref{d-3} and \eqref{d-5}, we have that
$$
\bar{\gamma}=\sup\left\{a>0 \  \big| \ \sup_{0 < s \le \bar\rho_0}s^{-a}p(s) < \infty \right\}.
$$
It follows from  \eqref{pcon1} and  \eqref{pcon2} that
\begin{align*}
0< K_2:=\inf_{0 < s\le 2\bar\rho_0} \frac{sp'(s)}{p(s)} \le \sup_{0 < s\le 2\bar\rho_0} \frac{sp'(s)}{p(s)}=:K_3<\infty.
 \end{align*}
Using \eqref{pcon1}, \eqref{d-0}, the L'Hopital rule and \eqref{pcon2},  one has that
\begin{align}
& \lim_{s\to 0+}\frac{sp''(s)}{p'(s)}= \lim_{s\to 0+}\frac{sp'(s)}{p(s)}-1 = \bar{\gamma} -1 \in \lt[3^{-1}, \infty\rt),\notag\\
&\sup_{0 < s\le 2\bar\rho_0}
\frac{s|p''(s)|}{p'(s)}=:K_4 \in (0, \infty). \notag
\end{align}
So, it holds that for $s\in (0, 2\bar\rho_0]$,
\begin{align}\label{21.3.15}
 K_2 p(s) \le  sp'(s) \le K_3 p(s) \ {\rm and} \
 s|p''(s)| \le K_4 p'(s).
\end{align}

In view of  \eqref{d-00}, we see that
\begin{align}
& i(0)=0, \label{21-11-2}\\
& i(s)=\int_0^s \tau^{-1}p'(\tau) d\tau
=s^{-1}p(s)+ \int_0^s \tau^{-2} p(\tau) d\tau,\notag
\end{align}
which gives, with the aid of \eqref{pcon2}, that for $s\in (0, \bar\rho_0]$,
\begin{align*}
&i(s)\le  s^{-1}p(s)+ \frac{3}{4} \int_0^s \tau^{-1} p'(\tau) d\tau=s^{-1}p(s)+ \frac{3}{4} i(s) .
\end{align*}
Then,  we obtain that for $s\in (0, \bar\rho_0]$,
\begin{align}\label{irho}
  p(s) \le  s i(s) \le 4 p(s),
\end{align}
which, together with  $p'>0$ and $i'>0$, implies that for
$s\in (0, 2\bar\rho_0]$,
\begin{align}\label{irho-1}
  p(s) \le  s i(s) \le K_5 p(s),
\end{align}
where $K_5=\max\{4,  \  2\bar\rho_0 i(2\bar\rho_0)/p(\bar\rho_0) \}$.
It follows from   \eqref{irho} and \eqref{d-3} that for $a<\bar\gamma-1$,
\begin{align*}
\lim_{s\to 0+} \frac{i(s)}{s^a}=  \lim_{s\to 0+} \frac{si(s)}{p(s)}\frac{p(s)}{s^{a+1}}\le  4\lim_{s\to 0+} \frac{p(s)}{s^{a+1}}=0,
\end{align*}
and
from   \eqref{irho} and \eqref{d-5} that for $a>\bar\gamma-1$,
\begin{align*}
\lim_{s\to 0+} \frac{s^a}{i(s)}=  \lim_{s\to 0+} \frac{s^{a+1}}{p(s)} \frac{p(s)}{si(s)}\le  \lim_{s\to 0+} \frac{s^{a+1}}{p(s)}=0.
\end{align*}
This  means that for $s\in (0, \bar\rho_0]$,
\begin{align}\label{April4}
 i(s)\le K_6  s^a   \  {\rm when} \   a<\bar\gamma-1,  \ \
s^a \le K_7  i(s)  \  {\rm when} \   a>\bar\gamma-1,
\end{align}
where  $K_6=\sup_{0<s\le \bar\rho_0} s^{-a}i(s)$ and  $K_7=\sup_{0<s\le \bar\rho_0} s^ai^{-1}(s)$ are finite positive constants.

We present some properties for the stationary solution $\bar\rho(x)$.
It follows from \eqref{pxphi} and $i(\bar\rho(\bar R))=i(0)=0$ that for $x\in I$,
\begin{align}\label{har-1}
K_{8}(\bar R -x) \le i(\bar\rho(x))=\int_x^{\bar R} s\phi(s)ds \le  K_{9} (\bar R -x),
\end{align}
where
$K_{8}= {M}\mathcal{G}/({2 \bar R^2})$ and $K_{9}={4\pi \mathcal{G}\bar\rho_0 \bar R}/{3}$.
Let
$$\tilde \rho(x,t)=\bar\rho(x) \left\{1+\vartheta \left(\frac{x^2}{r^2 r_x}
-1\right)\right\} \ {\rm with} \ \vartheta\in [0,1],$$
then it produces from \eqref{priori1} and the smallness of $\varepsilon_0$ that
\begin{align}\label{til-1}
(1/2)\bar\rho(x) \le \tilde \rho(x,t) \le (3/2) \bar\rho(x).
\end{align}
Notice that $|p(\tilde \rho )-p(\bar\rho)|=p'(\underline{\rho})|\tilde \rho-\bar\rho|$ with $\underline{\rho}=\vartheta_1 \tilde \rho + (1-\vartheta_1) \bar\rho$
for some constant $\vartheta_1\in(0,1)$, which implies, using \eqref{21.3.15}, $p'>0$ and
\eqref{priori1}, that
\begin{align*}
& |p(\tilde \rho )-p(\bar\rho)|
\le K_3  p(\underline{\rho})
|\underline{\rho}^{-1}(\tilde \rho-\bar\rho)|\\
&
\le  8 K_3 \epsilon_0 \max\{p(\tilde \rho ), \ p(\bar\rho)\}
\le 2^{-1}\max\{p(\tilde \rho ), \ p(\bar\rho)\}.
\end{align*}
Thus,
\begin{align}
 & 2^{-1}{p(\bar\rho)}\le  {p(\tilde \rho )}\le 2{p(\bar\rho)}, \label{til-2}\\
 &16^{-1}{i(\bar\rho)}\le  {i(\tilde \rho )}\le 4K_5 {i(\bar\rho)},\label{til-3}
\end{align}
where $\eqref{til-3}$ follows from \eqref{irho}, \eqref{irho-1}, \eqref{til-1} and \eqref{til-2}.

{\em The Hardy inequality}. Let $k>1$ be a given real number, $l$   a positive constant, and $f$   a function satisfying  $\int_0^{l} x^k(f^2+|f'|^2)dx <\infty$, then it holds that
\begin{equation}\label{hardy1}
\int_0^{l} x^{k -2} f^2 dx \leq \bar C\int_0^{l} x^{k}\lt( f^2+|f'|^2\rt)dx
\end{equation}
for a certain constant $\bar C$ depending only on $k$ and $l$, whose proof can be found in \cite{KMP}. In fact, \eqref{hardy1} is a general version of the standard Hardy inequality: $\int_0^\infty |x^{-1}f|^2 dx\le C \int_0^\infty |f'|^2 dx$ for some constant $C$. As a consequence of \eqref{hardy1},
we have the following estimates: for $l\in(0, \bar R)$,
\begin{equation}\label{hardy2}
\int_{\bar R-l}^{\bar R} (\bar R -x)^{k-2} f^2 dx\le \bar C \int_{\bar R-l}^{\bar R} (\bar R -x)^k\lt( f^2+|f'|^2\rt)dx<\infty.
\end{equation}

{\em Notation}. throughout the rest of this paper, $C$ will denote a positive constant which does not depend on the time $t$ or the data, and $C(\beta)$   a certain positive number depending only on quantity $\beta$. They are referred as universal and can change from one inequality to another one. We will adopt the notation $ a \les b$ to denote $a \le C b$,  where $C$ is the universal constant as defined
above. We will use $\int$ to denote $\int_I$.

\subsection{Lower-order estimates}\label{22.1.4-2}
This subsection consists of Lemmas \ref{lowt1}-\ref{lempoint}.

\begin{lemma}\label{lowt1}
It holds that for $t\in [0,T]$,
\begin{align}
&\mathcal{E}_0(t)+
(1+t)\sum_{i=1,2} \left(\mathcal{E}_i+\mathcal{D}_{i-1}\right)(t)
+ \int_0^t\mathcal{D}_0(s) ds\notag\\
& +\sum_{i=1,2}\int_0^t  (1+s)  \mathcal{D}_i(s) ds\les \sum_{0\le i\le 2}\mathcal{E}_i(0) +\sum_{i=0,1}\mathcal{D}_i(0),
\label{vvt}
\end{align}
where
\begin{align*}
&\mathcal{E}_0(t)=\int  \left[(r-x)^2+x^2 (r_x-1)^2  \right] dx,\\
&\mathcal{D}_0(t)=\int p(\bar\rho)\left[(r-x)^2
+x^2 (r_x-1)^2\right]dx,\\
&\mathcal{E}_i(t)=\int x^2\bar\rho \left|\partial_t^{i}r\right|^2 dx, \ \ \mathcal{D}_i(t) =\int   \lt(\left|\partial_t^{i}r\right|^2 +\left|x\partial_t^{i}r_x\right|^2\rt)dx, \  i= 1,2.
\end{align*}

\end{lemma}

{\em Proof}. To prove the lemma, we  use \eqref{pxphi} to
rewrite \eqref{fixedp-a}  as
\begin{align}\label{remaineq}
\bar\rho \frac{x^2}{r^2} v_t+\lt[p\lt(\frac{x^2 \bar\rho}{r^2 r_x}
\rt)\rt]_x-\frac{x^4}{r^4} \lt(p(\bar\rho)\rt)_x=\mathscr{B}_x+4\nu_1 \lt(\frac{v}{r}\rt)_x,
\end{align}
where $\mathscr{B}$ is defined by \eqref{mb}.
The proof consists of three steps.

{\em Step 1}. In this step, we prove that
\begin{equation}\label{l2t0}
(\mathcal{E}_1+\mathcal{D}_0)(t)
+\int_0^t  \mathcal{D}_1(s)  ds
\les (\mathcal{E}_1+\mathcal{D}_0)(0).
\end{equation}
We set
$
A(s)=\int_0^s \tau^{-2 }{p(\tau)}d\tau$
 for $  s >0$,
and integrate the product of \eqref{remaineq} and $r^2 v$ over $I$ to get
\begin{equation}\label{l2}
\frac{1}{2}\frac{d}{dt} \mathcal{E}_1(t)+
\frac{d}{dt}\int x^2 \eta(x,t) dx +\int D(x, t) dx =0,
\end{equation}
where
\begin{subequations}\label{10.14.2}\begin{align}
&\eta(x,t)= \bar\rho  A\biggl(\frac{x^2}{r^2}
\frac{\bar\rho}{r_x}\biggr) + p(\bar\rho)\biggl(\frac{x^2}{r^2}r_x-4\frac{x}{r}\biggr)
 - \bar\rho A(\bar\rho) +3  p(\bar\rho),  \\
& D(x,t)= \mathscr{B}(r^2 v)_x-4\nu_1\lt(\frac{v}{r}\rt)_xr^2 v\geq 3\sigma \lt(\frac{r^2}{r_x}v_x^2+ 2r_x v^2\rt)\label{double}
\end{align}\end{subequations}
with $\sigma=\min\{2\nu_1/3,\ \nu_2\}$. Indeed, \eqref{double} follows from the estimates obtained in \cite{LXZ2}.

It follows from the Taylor expansion that
\begin{equation}\label{h0}
\eta( x,t)=
  2^{-1} G(x,t)
+ \widetilde{G}(x,t),
\end{equation}
where
\begin{align}
& G=\bar\rho p'(\bar\rho)(r_x-1)^2+4\left(\bar\rho p'(\bar\rho)
-p(\bar\rho)\right)(r/x-1)^2\notag\\
& \ \ +4\left(\bar\rho p'(\bar\rho)-2p(\bar\rho)\right)(r/x-1)(r_x-1),
\label{G} \\
& | \widetilde{G} | \lesssim \left\{ p(\bar\rho)+ \bar\rho p'(\bar\rho) + \left( {\bar\rho}/{\tilde\rho} \right)^4\left| \tilde{\rho}^2 p''(\tilde{\rho}) - 4 \tilde{\rho} p'(\tilde{\rho} )    + 6p(\tilde{\rho})  \right| \right\}\notag\\
& \ \ \times
\left(\left|r_x-1\right|^3
+\left|\frac{r}{x}-1\right|^3 \right) \ {\rm with} \
\tilde \rho=\bar\rho \left\{1+\vartheta_2 \left(\frac{x^2}{r^2}
\frac{1}{r_x}-1\right)\right\} \notag
\end{align}
for some constant $\vartheta_2\in (0,1)$.
Simple calculations show that
\begin{align}\label{jo.1}
\frac{3G}{4p(\bar\rho)}=\left(r_x-\frac{r}{x}\right)^2
+\left(\frac{3\bar\rho p'(\bar\rho)}{4p(\bar\rho)} - 1\right) \left(
\left(r_x-\frac{r}{x}\right)+3\left(\frac{r}{x}-1\right) \right)^2,
\end{align}
which, together with \eqref{21.3.15}, means
\begin{align}
 \int x^2 G(x,t) dx \les  \mathcal{D}_0(t). \label{12.30-1}
\end{align}
In view of  \eqref{ibm} and \eqref{pcon2}, we see  that for $x\in [0, \bar R]$,
\begin{align}\label{basic1}
 \frac{3 \bar\rho(x) p'(\bar\rho(x))}{4p(\bar\rho(x))} \ge  1, \ {\rm and}  \  \lim_{x\to \bar R}\frac{3 \bar\rho(x) p'(\bar\rho(x))}{4p(\bar\rho(x))} =\frac{3}{4}\bar{\gamma}>1.
\end{align}
Thus, there exists a constant $\iota\in (0, \bar R/4]$ such that for $x\in [\bar R-2\iota, \bar R]$,
\begin{align}\label{basic2}
 \underline{c}=\frac{1}{2} \left(\frac{3}{4}\bar{\gamma}-1\right)\le \frac{3 \bar\rho(x) p'(\bar\rho(x))}{4p(\bar\rho(x))} -1 \le 3 \underline{c},
\end{align}
which means, with the help of \eqref{jo.1}, that for $x\in [\bar R-2\iota, \bar R]$,
\begin{align*}
\frac{3G}{4p(\bar\rho)}
\ge\left(r_x-\frac{r}{x}\right)^2
+\underline{c}\left(
\left(r_x-\frac{r}{x}\right)+3\left(\frac{r}{x}-1\right) \right)^2 .
\end{align*}
Thus, one has that for $x\in [\bar R-2\iota, \bar R]$,
\begin{subequations}\label{jo.2}\begin{align}
&p(\bar\rho)\left({r}/{x}-1\right)^2 \le
6^{-1}(1+\underline{c}^{-1})G ,\label{jo.2-a}\\
&p(\bar\rho)\left(r_x-1\right)^2 \le
6^{-1}\left(11+2\underline{c}^{-1}\right)G. \label{jo.2-b}
\end{align}\end{subequations}
Away from the vacuum boundary, it follows from \eqref{jo.1} and \eqref{basic1} that
$
3G/4 \ge   p(\bar\rho) (r_x-r/x)^2
$,
which, together with $p(\bar\rho(x))\ge p(\bar\rho(\bar R- \iota))>0$ on $[0, \bar R- \iota]$, implies
\begin{align}
&\frac{3}{4}\int \chi x^2 G dx \ge p(\bar\rho(\bar R- \iota)) \int \chi  (x r_x-r)^2 dx \notag \\
 = & p(\bar\rho(\bar R- \iota)) \lt\{\int \chi \lt[  (x r_x-x)^2 + 2(r-x)^2  \rt]dx
+\int   \chi' x (r-x)^2    dx\rt\}
\notag\\
\ge & p(\bar\rho(\bar R- \iota)) \int^{\bar R-2\iota}_0   \lt[  (x r_x-x)^2 + 2(r-x)^2  \rt]dx \notag\\
&-\frac{4}{\iota}\int_{\bar R-2\iota}^{\bar R- \iota}  p(\bar\rho)  x (r-x)^2    dx , \label{jo.3}
\end{align}
where $\chi=\chi(x)$ is a smooth cut-off function satisfying $\chi=1$ on $[0, \bar R-2\iota]$, $\chi=0$ on $[ \bar R-\iota, \bar R]$, and $-4/\iota\le \chi' \le 0$. It follows from \eqref{jo.2} and \eqref{jo.3} that
$
\mathcal{D}_0(t)\les \int x^2 G(x,t) dx  $. This, together with \eqref{12.30-1}, gives
\begin{align}
 \mathcal{D}_0(t)\les \int x^2 G(x,t) dx \les \mathcal{D}_0(t).\label{h1}
\end{align}
We  use \eqref{priori1},
\eqref{21.3.15} and \eqref{til-1} to obtain
\begin{align*}
| \widetilde{G} |\lesssim \epsilon_0 p(\bar\rho) \left( 1+ \frac{p(\tilde{\rho})}{p(\bar\rho)}\right) \left(\left(r_x-1\right)^2 + \left(\frac{r}{x}-1\right)^2\right),
 \end{align*}
which means, with the aid of   \eqref{til-2}, that
\begin{align}\label{h3}
\int x^2 | \widetilde{G}(x,t) | dx \lesssim \epsilon_0 \mathcal{D}_0(t) .
 \end{align}

We integrate the product of $(1+t)^{\ell}$ $(\ell\ge 0)$ and \eqref{l2} over $[0, t]$, and use \eqref{double},
\eqref{h0},   \eqref{h1} and \eqref{h3}  to  achieve
\begin{align}\label{l2tl}
&(1+t)^\ell(\mathcal{E}_1+\mathcal{D}_0)(t)
+\int_0^t (1+s)^\ell \mathcal{D}_1(s)  ds\notag\\
&\les (\mathcal{E}_1+\mathcal{D}_0)(0)
+ \ell \int_0^t (1+s)^{\ell-1} (\mathcal{D}_0+\mathcal{D}_1)(s) ds,
\end{align}
because of $\mathcal{E}_1\les \mathcal{D}_1$ which is due to \eqref{yo.1}.
Letting $\ell=0$ in \eqref{l2tl} proves \eqref{l2t0}.

{\em Step 2}. In this step, we prove that
\begin{equation}\label{pr2v}
\begin{aligned}
&\mathcal{E}_0(t) + (1+t)(\mathcal{E}_1+ \mathcal{D}_0)(t) +
\int_0^t \left[\mathcal{D}_0(s) \right.\\
&\left.+ (1+s) \mathcal{D}_1(s) \right] ds
\les (\mathcal{E}_0+\mathcal{E}_1+ \mathcal{D}_0)(0).
\end{aligned}
\end{equation}
It follows from \eqref{l2tl} with
 $\ell=1$ and \eqref{l2t0} that
\begin{align}
(1+t)(\mathcal{E}_1+\mathcal{D}_0)(t)
+ \int_0^t(1+s) \mathcal{D}_1(s)  ds\notag\\
\les (\mathcal{E}_1+\mathcal{D}_0)(0) + \int_0^t \mathcal{D}_0(s) ds.\label{l2t1}
\end{align}
To bound the last term in \eqref{l2t1},
 we
integrate the product of \eqref{remaineq} and $r^3- x^3$  over $I$ to get
\begin{align}
&\frac{d}{dt}\int \left\{ x^2\eta_0(x,t)+x^3\bar\rho v\left({r}/{x}- {x^2}/{r^2}\right)\right\}dx
\notag\\
& +\int D_0(x,t) dx
=\int x^2\bar\rho v^2\biggl(1+2{x^3}/{r^3}\biggr)dx,\label{r3x3}
\end{align}
where
\begin{subequations}\label{10.14.1}\begin{align}
&\eta_0= 4\nu_1 \lt[\ln\lt(\frac{r}{xr_x}\rt)
+\frac{xr_x}{r}-1\rt]
+3\nu_2\lt[\frac{r^2r_x}{x^2}-\ln\lt(\frac{r^2 r_x}{x^2}\rt)-1\rt], \\
&D_0=p(\bar\rho)\biggl[\frac{x^4}{r^4}(r^3-x^3)\biggr]_x
-
p\biggl(\frac{x^2\bar\rho}{r^2r_x}
\biggr)(r^3-x^3)_x.
\end{align}\end{subequations}
Clearly, we can use the estimates obtained in \cite{LXZ2} to show
\begin{align}
   \eta_0  \les |r/x-1|^2+|r_x-1|^2  \les  \sigma^{-1}  \eta_0  , \label{eta0}
\end{align}
where $\sigma=\min\{2\nu_1/3,\ \nu_2\}$.
It follows from the Taylor expansion that
$
D_0=3x^2 G+ x^2 \widetilde{G}_0$,
where $G$ is defined by \eqref{G}, and
$$
 |\widetilde{G}_0| \lesssim
\left( p(\bar\rho)+\bar\rho p'(\bar\rho) + \left( {\bar\rho}/{\hat\rho} \right)^2 \hat{\rho}^2 | p''(\hat{\rho})| \right) \left(\left|r_x-1\right|^3 +\left|{r}/{x}-1\right|^3 \right)
$$
with $\hat \rho=\bar\rho \left\{1+\vartheta_3 \left({x^2}/({r^2r_x})
-1\right)\right\}$
for some constant $\vartheta_3\in (0,1)$.
In a similar way to deriving \eqref{h1} and \eqref{h3}, one has
\begin{align}
\mathcal{D}_0(t)\les \int x^2 G(x,t) dx \le \int D_0(x,t) dx. \label{r3x3est}
\end{align}
We integrate \eqref{r3x3} over $[0,t]$, and use   \eqref{l2t0},   \eqref{eta0}, \eqref{r3x3est} and the Cauchy inequality to obtain
\begin{align}\label{l2r3x3}
\mathcal{E}_0(t)+\int_0^t \mathcal{D}_0(s) ds \les
(\mathcal{E}_0+ \mathcal{E}_1+\mathcal{D}_0)(0),
\end{align}
where we have used the following estimates:
\begin{align*}
&\int x^2\bar\rho v^2\biggl(1+2{x^3}/{r^3}\biggr)dx \les \mathcal{E}_1(t) \les  \mathcal{D}_1(t),\\
&\int \left| x^3\bar\rho v\left({r}/{x}- {x^2}/{r^2}\right) \right| dx \les \int \left| x^2 \bar\rho v (r-x)\right| dx\\
&\les \omega \int  (r-x)^2 dx
+\omega^{-1} \int x^4 \bar\rho^{2} v^2 dx \les \omega \mathcal{E}_0(t) +\omega^{-1} \mathcal{E}_1(t)
\end{align*}
for any $\omega>0$. As a consequence of \eqref{l2t1} and \eqref{l2r3x3}, we prove \eqref{pr2v}.

{\em Step 3}. In this step, we prove
\begin{align}\label{l2tderi}
(1+t)\left(\mathcal{E}_2+\mathcal{D}_1\right)(t)
+ \int_0^t(1+s)\mathcal{D}_2(s) ds
\les \sum_{0\le i\le 2}\mathcal{E}_i(0) +\sum_{i=0,1}\mathcal{D}_i(0).
\end{align}
 Setting $\varrho=\bar\rho x^2/(r^2 r_x)$, and
integrating the product of  $\partial_t(r^2\eqref{remaineq})$ and  $v_t$ over $I$, we have
\bee\label{phivt}
\frac{1}{2}\frac{d}{dt}\mathcal{E}_2(t)
+\frac{d}{dt}\int \Phi(x, t) dx +\int J_1(x,t) dx=\int J_2(x,t) dx,
\ene
where
\begin{align*}
&\Phi=\left(2\varrho p'\left(\varrho\right)-p \left(\varrho\right)\right)r_x v^2+2\left(\varrho p'\left(\varrho\right)-p \left(\varrho\right)\right)rvv_x
\\
&\quad +\frac{r^2}{2r_x}\varrho p'\left(\varrho\right)v_x^2 -p(\bar\rho)\biggl[\biggl(4\frac{x^3}{r^3}-3\frac{x^4}{r^4}r_x\biggr)v^2+2\frac{x^4}{r^3}v v_x\biggr],\\
&J_1=\lt[\mathscr{B}_t(r^2 v_t)_x- 4\nu_1r^2v_t\lt(\frac{v}{r}\rt)_{xt}\rt]
 +2 \left[ \mathscr{B}(rvv_t)_x-4\nu_1\lt(\frac{v}{r}\rt)_xrvv_t \right],\\
&J_2=\left[\lt(2\varrho p'\left(\varrho\right)-p \left(\varrho\right)\rt)r_x\right]_t v^2+2\left[\lt(\varrho p'\left(\varrho\right)-p \left(\varrho\right)\right)r\rt]_tvv_x
\\
&\quad +\lt[\frac{r^2}{2r_x}\varrho p'\left(\varrho\right)\rt]_tv_x^2 -p(\bar\rho)
\lt[\lt(4\frac{x^3}{r^3}-3\frac{x^4}{r^4}r_x\rt)_tv^2
+2\lt(\frac{x^4}{r^3}\rt)_t v v_x\rt].
\end{align*}

Clearly, it holds that
\begin{subequations}\label{22.1.1-1}\begin{align}
&J_1  \ge  2   \sigma   \lt[ {r^2} v_{tx}^2/r_x+2r_x v_t^2\rt]
-C \sigma^{-1}  \lt(x^2 v_x^2 + v^2\rt), \label{21.10.13-1}\\
&| J_2 | \les  x^2 v_x^2 + v^2,  \label{21.10.13-2}
\end{align}\end{subequations}
where $\sigma=\min\{2\nu_1/3,\ \nu_2\}$. Indeed, \eqref{21.10.13-1} follows from the estimates obtained in \cite{LXZ2}; and \eqref{21.10.13-2} from \eqref{pcon1-a}, \eqref{ibm},  \eqref{21apri} and \eqref{21.3.15}.
Note that
\begin{equation}\label{21.10.13-3}
\Phi(x, t)= 2^{-1}G_1(x,t)+\widetilde \Phi(x,t),
\end{equation}
where
\begin{align*}
& G_1=\bar\rho p'(\bar\rho) x^2v_x^2+ 4\lt(\bar\rho p'(\bar\rho)-p(\bar\rho)\rt) v^2+4\lt (\bar\rho p'(\bar\rho)-2p(\bar\rho)
\rt)vxv_x,
\notag\\
&  \widetilde \Phi=\left[\lt(2\varrho p'\left(\varrho\right)-p \left(\varrho\right)\rt)r_x-\lt(2\bar\rho p'(\bar\rho)-p(\bar\rho)\rt) \right]  v^2\\
&\quad +2\left[\lt(\varrho p'\left(\varrho\right)-p \left(\varrho\right)\right)r -\lt(\bar\rho p'\left(\bar\rho\right)-p \left(\bar\rho\right)\right)x \rt] v v_x
\\
&\quad +\frac{1}{2}\lt[\frac{r^2}{r_x}\varrho p'\left(\varrho\right)-x^2 \bar\rho p'\left(\bar\rho\right)\rt] v_x^2 \\
& \quad -p(\bar\rho)
\lt[\lt(4\frac{x^3}{r^3}-3\frac{x^4}{r^4}r_x-1\rt) v^2
+2\lt(\frac{x^4}{r^3}-x\rt)  v v_x\rt].
\end{align*}
Then, we use a similar way to deriving \eqref{h1} and \eqref{h3} to get
\begin{subequations}\label{22.1.1-2}\begin{align}
& \int G_1(x, t)dx \les \int  p(\bar\rho) (v^2 +x^2v_x^2)dx \les \int G_1(x, t)dx ,\label{G-1}\\
&\int |\widetilde\Phi(x,t)| dx \les \epsilon_0  \int  p(\bar\rho) (v^2 +x^2v_x^2)dx. \label{tphi}
\end{align}\end{subequations}
We integrate \eqref{phivt} over $[0,t]$, and use \eqref{22.1.1-1}-\eqref{22.1.1-2} and \eqref{pr2v}   to obtain
\begin{align}
&\mathcal{E}_2(t)+\int  p(\bar\rho) (v^2 +x^2v_x^2)dx
+\int_0^t \mathcal{D}_2(s)ds  \les \sum_{0\le i\le 2}\mathcal{E}_i(0) +\sum_{i=0,1}\mathcal{D}_i(0). \label{21-10-13}
\end{align}
Similarly, we integrate the product of $1+t$ and \eqref{phivt} over $[0,t]$ to achieve
\begin{align}
&(1+t)\mathcal{E}_2(t)+(1+t)\int  p(\bar\rho) (v^2 +x^2v_x^2)dx
+\int_0^t (1+s) \mathcal{D}_2(s)ds \notag\\
& \les \sum_{0\le i\le 2}\mathcal{E}_i(0) +\sum_{i=0,1}\mathcal{D}_i(0), \label{21-10-13-1}
\end{align}
where \eqref{21-10-13} has been used. Note that for any function $f$ and positive constant $a$,
\begin{align}\label{12.13.a}
& (1+t)^a\int f^2(x,t) dx \le  \int f^2(x,0) dx
+\int_0^t (1+s)^a \int (f^2 +f_s^2)(x,s)dxds\notag\\
&\qquad
+a\int_0^t (1+s)^{a-1} \int f^2(x,s)dxds,
\end{align}
then \eqref{l2tderi} follows from \eqref{pr2v} and
\eqref{21-10-13-1}.
\hfill $\Box$

\begin{lemma}\label{l2palpha}
Let $\Upsilon$ be defined by
\eqref{22.1.11-1}. Then it holds that for
$\theta\in \lt(0, 1-5/(4\bar\gamma)\rt]$ and $t\in [0, T]$,
\begin{subequations}\label{12-21}\begin{align}
&\mathscr{E}_0(t)+(1+t)^\zeta\mathscr{D}_0(t)
+(1+t)^{2\zeta}\left(
\mathcal{D}_0+\mathcal{D}_1+\mathcal{E}_2\right)(t)
+(1+t)^{2\zeta-1}\mathcal{E}_0(t)\notag\\
&\quad+\int_0^t\lt[\mathscr{D}_0(s)
+(1+s)^\zeta\mathscr{D}_1(s)
+(1+s)^{2\zeta}(\mathcal{D}_1+\mathcal{D}_2)
(s)\rt. \notag\\
&\quad\lt. +
(1+s)^{2\zeta-1}\mathcal{D}_0(s)\rt]ds
 \les  C(\theta^{-1}) \lt(\mathscr{E}_0 +\mathcal{D}_1+\mathcal{E}_2\rt)(0),
 \label{l2palphaest}\\
& \lt( \Upsilon(1+t)^{3- {3}/{\bar\gamma}-2\theta} + (1-\Upsilon) (1+t)^{2\zeta} \rt) \int (r(x, t)-x)^2 dx\notag\\
&\quad\les C(\theta^{-1}) \lt(\mathscr{E}_0 +\mathcal{D}_1+\mathcal{E}_2\rt)(0), \label{rxt00}
\end{align}\end{subequations}
where  $\alpha=1-\theta$, $\zeta=1- {1}/({2\bar\gamma})
- \theta/2$,
\begin{align*}
&\mathscr{E}_0(t)=\int  i^{-\alpha}(\bar\rho)\left[(r-x)^2+x^2(r_x-1)^2
\right]dx,\\
&\mathscr{D}_0(t)=\int i^{-\alpha}(\bar\rho) p(\bar\rho)\left[(r-x)^2+x^2(r_x-1)^2
\right]dx,\\
&\mathscr{D}_1(t)=\int  i^{-\alpha}(\bar\rho)(v^2+x^2v_x^2) dx.
\end{align*}

\end{lemma}

{\em Proof}. The proof consists of four steps.

{\em Step 1}.  In this step, we prove that for any $\omega>0$,
\begin{align}
&\mathscr{E}_0(t)+\int_0^t \mathscr{D}_0(s) ds  \les \mathscr{E}_0(0)
+\lt(1+\omega^{-1}\rt)\sum_{0\le i\le 2}\int_0^t \mathcal{D}_i(s)ds
\notag \\
&+ \omega \int_0^t  \lt\{\mathscr{D}_0(s)+(1+s)^{\zeta} \mathscr{D}_1(s)
+(1+s)^{2\zeta} \mathcal{D}_1(s)\rt\} ds\notag\\
& + \int_0^t  \lt\{ \omega^{1-2q_1} (1+s)^{- q_1 \zeta}+\omega^{-1} \theta^{-1}(1+s)^{-2\zeta} \rt\}\mathscr{E}_0(s)ds,
\label{imt0}
\end{align}
where $q_1=2\bar\gamma/(2\bar\gamma-1- 2\bar\gamma \theta)$ satisfies
 $q_1 \zeta>1$.

We integrate the product of  \eqref{remaineq} with $\int_0^x i^{-\alpha}(\bar\rho)(r^3-y^3)_y \ dy$  over $I$ to get
\begin{equation}
\frac{d}{dt}\int i^{-\alpha}(\bar\rho) x^2 \eta_0 dx + \int i^{-\alpha}(\bar\rho) D_0 dx =\sum_{1\le i\le 3} L_i, \label{22.1.2-1}
\end{equation}
where $\eta_0$ and $D_0$ are defined by \eqref{10.14.1}, and
\begin{align*}
&L_1=-\int \bar\rho \frac{x^2}{r^2} v_t\int_0^x i^{-\alpha}(\bar\rho)(r^3-y^3)_y
dydx,\\
&L_2=\int p(\bar\rho)\biggl(\frac{x^4}{r^4}\biggr)_x\left[i^{-\alpha}(\bar\rho)(r^3-x^3)-\int_0^x i^{-\alpha}(\bar\rho)(r^3-y^3)_ydy \right]dx,\\
&L_3=4\nu_1 \int \lt(\frac{v}{r}\rt)_x\lt[\int_0^x i^{-\alpha}(\bar\rho)(r^3-y^3)_y dy -i^{-\alpha}(\bar\rho)(r^3-x^3)\rt]dx.
\end{align*}
In a similar way to deriving \eqref{r3x3est}, one has
$$\mathscr{D}_0(t)\les \int i^{-\alpha}(\bar\rho) x^2 G(x,t) dx \le  \int i^{-\alpha}(\bar\rho) D_0(x,t) dx. $$
Integrating \eqref{22.1.2-1} over $[0,t]$, and using \eqref{eta0} and the inequality above, we have
\begin{equation}\label{r3x3alpha}
\mathscr{E}_0(t)+\int_0^t \mathscr{D}_0(s) ds  \les \mathscr{E}_0(0)+\sum_{1\le i\le 3} \int_0^t | L_i| ds.
\end{equation}

For $L_1$, it follows from the Cauchy and H$\ddot{o}$lder inequalities that for any $\omega>0$,
\begin{align*}
|L_1|\les & \omega^{-1}\int v_t^2 dx+\omega \int\bar\rho^2  \left|\int_0^x i^{-\alpha}(\bar\rho)y\lt(|r-y|+y|r_y-1|\rt)dy\right|^2dx\notag\\
\les &  \omega^{-1} \mathcal{D}_2(t)+\omega \mathscr{D}_0(t) H_1,
\end{align*}
where
$H_1=\int \bar\rho^2 \int_0^xp^{-1}(\bar\rho) i^{-\alpha}(\bar\rho)y^2dy dx$.
Due to the L'Hopital rule, \eqref{pxphi} and \eqref{irho}, we have
$$\lim_{x\rightarrow \bar R}\frac{\int_0^xp^{-1}(\bar\rho) i^{-\alpha}(\bar\rho)dy}{p^{-1}(\bar\rho)i^{1-\alpha}(\bar\rho)}
=\lim_{x\rightarrow \bar R}\frac{1}{\lt(\bar\rho i(\bar\rho)p^{-1}(\bar\rho)+\alpha-1\rt)x\phi}\le \frac{{\bar R}^2}{\alpha M},
$$
which, together with \eqref{irho}, \eqref{April4} and \eqref{har-1}, means
\begin{align}\label{21.11.2-6}
H_1\les \int \bar\rho^2  p^{-1}(\bar\rho)i^{1-\alpha}(\bar\rho) dx \les
  \int \bar\rho  i^{-\alpha}(\bar\rho) dx
  \les \int    i^{\frac{8}{9(\bar\gamma-1)}-\alpha}(\bar\rho) dx\les 1.
\end{align}
So, we obtain that for any $\omega>0$,
\begin{align}\label{L1}
|L_1|\les \omega \mathscr{D}_0(t) +\omega^{-1} \mathcal{D}_2(t).
\end{align}

For $L_2$,
we can rewrite $L_2$ as $L_2=L_{21}+L_{22}$,
where
\begin{align}
L_{21}&=\int_0^{\bar R/2} p(\bar\rho)\lt(\frac{x^4}{r^4}\rt)_x\int_0^x
\left(i^{-\alpha}(\bar\rho)\right)_y(r^3-y^3)dy dx,\notag\\
L_{22}&=\int_{\bar R/2}^{\bar R} p(\bar\rho)\biggl(\frac{x^4}{r^4}\biggr)_x\left[i^{-\alpha}(\bar\rho)(r^3-x^3)-\int_0^x i^{-\alpha}(\bar\rho)(r^3-y^3)_ydy \right]dx.\notag
\end{align}
Clearly, it holds that
\begin{align}\label{L21}
|L_{21}|\les  \mathcal{D}_0(t).
\end{align}
(The derivation can be found in  \cite{LZeng}.)
It follows from the Cauchy and H$\ddot{o}$lder inequalities that for any $\omega>0$,
\begin{align}
&|L_{22}|\les\int_{\bar R/2}^{\bar R} p(\bar\rho)i^{-\alpha}(\bar\rho)(x|r_x-1|+|r-x|)|r-x|dx\notag\\
& +\int_{\bar R/2}^{\bar R} p(\bar\rho)(x|r_x-1|+|r-x|)
\int_0^xi^{-\alpha}
(\bar\rho)y(y|r_y-1|+|r-y|)dydx\notag\\
& \les 2 \omega \mathscr{D}_0(t)+\omega^{-1}H_2
 +\omega^{-1} \mathcal{D}_0(t) H_3, \label{L22}
\end{align}
where
\begin{align}
&H_2=\int_{\bar R/2}^{\bar R}p(\bar\rho)i^{-\alpha}(\bar\rho)(r-x)^2 dx, \notag\\
&H_3= \int_{\bar R/2}^{\bar R}p(\bar\rho)i^{\alpha}(\bar\rho)
\int_0^xp^{-1}(\bar\rho)i^{-2\alpha}
(\bar\rho)y^2dydx \les 1. \label{21.11.4}
\end{align}
Indeed, the bound for $H_3$ can be obtained by the same way as that for $H_1$ in \eqref{21.11.2-6}. 
In view of \eqref{April4}, we see that
$\bar\rho^{9(\bar\ga-1)/8}\les i(\bar\rho)\les \bar\rho^{7(\bar\ga-1)/8}$, which, together with \eqref{irho}, \eqref{har-1} and \eqref{hardy2}, gives

\begin{align}
&H_2\les \int_{\bar R/2}^{\bar R} i^{1-\alpha+\frac{8}{9(\bar\gamma-1)}}(\bar\rho)(r-x)^2 dx\notag \\
& \les
\int_{\bar R/2}^{\bar R} i^{3-\alpha+\frac{8}{9(\bar\gamma-1)}}(\bar\rho)
\left[(r-x)^2+  (r_x-1)^2  \right]dx\notag\\
& \les
\int_{\bar R/2}^{\bar R} \bar\rho^{\frac{7}{8}(\bar\gamma-1)(2-\alpha)-\frac{2}{9}}
p(\bar\rho)\left[(r-x)^2+  x^2 (r_x-1)^2  \right]dx \les \mathcal{D}_0(t).\notag
\end{align}
This gives, with the help of \eqref{L21}-\eqref{21.11.4}, that for any $\omega>0$,
\begin{equation}\label{L2}
|L_2|\les \omega\mathscr{D}_0(t)+\lt(1+\omega^{-1}\rt)\mathcal{D}_0(t).
\end{equation}

For $L_3$,
we can rewrite $L_3$ as
$L_3=-4\nu_1(L_{31}+L_{32})$, where
\begin{align}
&L_{31}=  \int_0^{\bar R/2 }\lt(\frac{v}{r}\rt)_x  \int_0^x (i^{-\alpha}( \bar\rho))_y(r^3-y^3)dy dx,\notag\\
& L_{32}=\int_{\bar R/2}^{\bar R} \lt(\frac{v}{r}\rt)_x\lt[ i^{-\alpha}(\bar \rho)(r^3-x^3)- \int_0^x i^{-\alpha}(\bar\rho)(r^3-y^3)_y dy \rt]dx.\notag
\end{align}
It is easy to show that
\begin{align}
|L_{31}|\les \lt(\mathcal{D}_1+\mathcal{D}_0\rt)(t), \ \
|L_{32}|\les L_{32}^I+L_{32}^{II},
\label{L31}
\end{align}
where
\begin{align*}
& L_{32}^I=\int_{\bar R/2}^{\bar R} i^{-\alpha}(\bar\rho) (x|v_x|+|v|)\left| r-x\right|dx ,\\
& L_{32}^{II}=\int_{\bar R/2}^{\bar R}(x|v_x|+|v|)\int_0^x i^{-\alpha}(\bar\rho)y\lt(y|r_y-1|+|r-y|\rt)dy dx.
\end{align*}
It follows from the Cauchy inequality, \eqref{har-1}, \eqref{hardy2} and $\alpha<1$ that for any $\omega>0$,
\begin{align}\label{L321}
L_{32}^I & \les  \omega (1+t)^{\zeta} \mathscr{D}_1(t)
+\omega^{-1} (1+t)^{-\zeta}
\int_{\bar R/2}^{\bar R} i^{ -\alpha}(\bar\rho) ( r-x)^2 dx\notag\\
& \les \omega  (1+t)^{\zeta} \mathscr{D}_1(t)
+\omega^{-1} (1+t)^{-\zeta}
\widetilde{L}_{32}^I ,
\end{align}
where
$$\widetilde{L}_{32}^I =\int_{\bar R/2}^{\bar R} i^{2-\alpha}(\bar\rho) \left[(r-x)^2+x^2(r_x-1)^2
\right]dx.$$
We choose   constants $q_1$, $q_2$ and $\delta_3$ satisfying
$$ \frac{1}{q_1}=1-\frac{1}{2\bar\gamma}-\theta, \ \  \frac{1}{q_2}=1-\frac{1}{q_1}
 \ \ {\rm and}  \  \ \delta_3=\frac{1}{2q_2-1},$$
then
$
 q_1, q_2>1$,  $q_1 \zeta>1$ and  $0<\delta_3<\bar\gamma-1$.
In view of \eqref{irho} and \eqref{April4}, we see that
$$i^{2q_2}(\bar\rho)\le 4 p(\bar\rho) \bar\rho^{-1} i^{2q_2-1}(\bar\rho) \les p(\bar\rho)  \bar\rho^{\delta_3 (2q_2-1)-1} =p(\bar\rho),$$
which, together with the H\"{o}lder inequality, implies
\begin{align}
\widetilde{L}_{32}^I \les\mathscr{E}_0^{1/q_1}
\lt(\int_{\bar R/2}^{\bar R} i^{2q_2-\alpha}(\bar\rho) \left[(r-x)^2+x^2(r_x-1)^2
\right]dx\rt)^{1/q_2}
\les \mathscr{E}_0^{1/q_1} \mathscr{D}_0^{1/q_2}.
\notag
\end{align}
Substitute this into \eqref{L321} and use the Young inequality  to get
\begin{align}\label{L321-21}
&L_{32}^I \les
 \omega \lt\{ \mathscr{D}_0(t)+(1+t)^{\zeta} \mathscr{D}_1(t) \rt\}
+ \omega^{1-2q_1} (1+t)^{- q_1 \zeta}\mathscr{E}_0(t) .
\end{align}
It follows from the Cauchy and H\"{o}lder inequalities that for any $\omega>0$,
\begin{align*}
L_{32}^{II} &\les \int_{\bar R/2}^{\bar R}(x|v_x|+|v|)\lt( \theta^{-1} \mathscr{E}_0(t) \rt)^{1/2}dx \\
& \les \omega (1+t)^{2\zeta} \mathcal{D}_1(t) + \omega^{-1} \theta^{-1}(1+t)^{-2\zeta} \mathscr{E}_0(t),
\end{align*}
which gives, with the help of \eqref{L31} and \eqref{L321-21}, that for any $\omega>0$,
\begin{align}\label{L3-}
|L_3|\les& (\mathcal{D}_0+ \mathcal{D}_1)(t)+
\omega \lt\{\mathscr{D}_0(t)+(1+t)^{\zeta} \mathscr{D}_1(t)
+(1+t)^{2\zeta} \mathcal{D}_1(t)\rt\}\notag\\
&+ \lt\{ \omega^{1-2q_1} (1+t)^{- q_1 \zeta}+\omega^{-1} \theta^{-1}(1+t)^{-2\zeta} \rt\}\mathscr{E}_0(t) .
\end{align}

Therefore, \eqref{imt0} follows from \eqref{r3x3alpha}, \eqref{L1}, \eqref{L2} and\eqref{L3-}.

{\em Step 2}.
In this step, we prove that for any $\omega>0$,
\begin{align}
&(1+t)^{\zeta}\mathscr{D}_0(t)
+\int_0^t(1+s)^\zeta \mathscr{D}_1(s)ds  \les  \mathscr{D}_0(0) +\int_0^t  \mathscr{D}_0(s) ds\notag\\
&+ (1+\omega^{-1}) \int_0^t (1+s) \lt(
 \theta^{-1} \mathcal{D}_1+\mathcal{D}_2 \rt)(s) ds\notag\\
& + \omega \int_0^t \lt\{ (1+s)^\zeta \mathscr{D}_1(s) + (1+t)^{2\zeta-1}\mathcal{D}_0(s) \rt\} ds .
\label{tzeta0}
\end{align}
We integrate the product of \eqref{remaineq} and $\int_0^x i^{-\alpha}(\bar\rho)(r^2 v)_y dy$  over $I$ to obtain
$$\frac{d}{dt}\int x^2 i^{-\alpha}(\bar\rho) \eta dx + \int i^{-\alpha}(\bar\rho) D dx = \sum_{1\le i\le 3} \mathcal{L}_i,$$
where $\eta$ and $D$ are given by \eqref{10.14.2}, and
\begin{align*}
&\mathcal{L}_1=-\int \bar\rho \frac{x^2}{r^2}v_t \int_0^x i^{-\alpha}(\bar\rho)(r^2 v)_ydydx,\\
&\mathcal{L}_2=\int p(\bar\rho)\lt(\frac{x^4}{r^4}\rt)_x\lt[ i^{-\alpha}(\bar\rho)r^2 v- \int_0^x i^{-\alpha}(\bar\rho) (r^2 v)_y dy\rt]dx,\\
&\mathcal{L}_3=4\nu_1\int  \left(\frac{v}{r}\right)_x\left[\int_0^x i^{-\alpha}(\bar\rho)(r^2 v)_ydy-i^{-\alpha}(\bar\rho)r^2 v\right]dx.
\end{align*}
We integrate the product of $(1+t)^\zeta$ and the equation above over $[0,t]$, and use the similar way to deriving \eqref{l2tl} to get
\begin{align}\label{zetakkk}
&(1+t)^{\zeta}\mathscr{D}_0(t)
+\int_0^t(1+s)^\zeta \mathscr{D}_1(s)ds\notag\\
& \les  \mathscr{D}_0(0)+ \sum_{i=1}^3\int_0^t (1+s)^\zeta |\mathcal{L}_i| ds+\int_0^t \mathscr{D}_0(s) ds.
\end{align}

The estimates for $\mathcal{L}_i$ $(i=1,2,3)$ are similar to those for $L_i$ $(i=1,2,3)$ obtained in the first step.
For $\mathcal{L}_1$, it follows from the Cauchy and  H\"{o}lder inequalities, and  \eqref{21.11.2-6} that for any $\omega>0$,
\begin{align}\label{k1}
|\mathcal{L}_1|  \les\omega \mathscr{D}_1(t)+\omega^{-1} \mathcal{D}_2(t),
\end{align}
due to
$\int \bar\rho^2 \int_0^x  i^{-\alpha}(\bar\rho) y^2 dy dx \les H_1 \les 1$.
For $\mathcal{L}_2$, note that
\begin{align*}
&\mathcal{L}_2=\int_{\bar R/2}^{\bar R} p(\bar\rho)\lt(\frac{x^4}{r^4}\rt)_x\lt[ i^{-\alpha}(\bar\rho)r^2 v- \int_0^x i^{-\alpha}(\bar\rho) (r^2 v)_y dy\rt]dx\\
&+\int_0^{\bar R/2} p(\bar\rho)\lt(\frac{x^4}{r^4}\rt)_x\int_0^x \lt(i^{-\alpha}(\bar\rho) \rt)_y r^2 v dydx=\mathcal{L}_{22}+\mathcal{L}_{21},
\end{align*}
and that for any $\omega>0$,
$
|\mathcal{L}_{21}|
\les\omega (1+t)^{\zeta-1}\mathcal{D}_0(t) +
\omega^{-1} (1+t)^{1-\zeta}\mathcal{D}_1(t)
$
and
\begin{align*}
&|\mathcal{L}_{22}|\les \int_{\bar R/2}^{\bar R} p(\bar\rho)i^{-\alpha}(\bar\rho)(|r-x|+x|r_x-1|)|v| dx\notag\\
&+\int_{\bar R/2}^{\bar R} p(\bar\rho)(|r-x|+x|r_x-1|)\int_0^x i^{-\alpha}(\bar\rho)y(|v|+|yv_y|)dydx\notag\\
& \les 2\omega (1+t)^{\zeta-1}\mathcal{D}_0(t)+
\omega^{-1} (1+t)^{1-\zeta} \lt\{ \int_{\bar R/2}^{\bar R} p(\bar\rho) i^{-2\alpha}(\bar\rho) v^2 dx\rt.\\
& \lt.   +  \mathcal{D}_1(t)
\int_{\bar R/2}^{\bar R} p(\bar\rho) \int_0^x i^{-2\alpha}(\bar\rho)y^2 dy dx \rt\} \\
& \les   2\omega (1+t)^{\zeta-1}\mathcal{D}_0(t) +
\omega^{-1} (1+t)^{1-\zeta}\mathcal{D}_1(t).
\end{align*}
Then,
 we have that for any $\omega>0$,
\begin{align}\label{k2}
|\mathcal{L}_{2}|\les  \omega (1+t)^{\zeta-1}\mathcal{D}_0(t) +
\omega^{-1} (1+t)^{1-\zeta}\mathcal{D}_1(t).
\end{align}
For $\mathcal{L}_3$, notice that
\begin{align*}
&-\frac{1}{4\nu_1}\mathcal{L}_3=\int_{\bar R/2}^{\bar R}  \left(\frac{v}{r}\right)_x\left[i^{-\alpha}(\bar\rho)r^2 v-\int_0^x i^{-\alpha}(\bar\rho)(r^2 v)_ydy\right]dx \\
&+\int^{\bar R/2}_0  \left(\frac{v}{r}\right)_x \int_0^x \lt(i^{-\alpha}(\bar\rho) \rt)_y r^2 v dydx= \mathcal{L}_{32} +\mathcal{L}_{31},
\end{align*}
and
$
|\mathcal{L}_{31}| \les \mathcal{D}_1(t)$,
$|\mathcal{L}_{32}| \les \mathcal{L}_{32}^{I}+ \mathcal{L}_{32}^{II}$,
where
\begin{align*}
\mathcal{L}_{32}^{I} &\les \int_{\bar R/2}^{\bar R}  i^{-\alpha}(\bar\rho) (x|v_x|+|v|)\left|v\right|dx\\
& \les \omega \mathscr{D}_1(t) + \omega^{-1} \int _{\bar R/2}^{\bar R}  i^{-\alpha}(\bar\rho) v^2 dx
\les \omega \mathscr{D}_1(t) \\
&\quad + \omega^{-1} \int _{\bar R/2}^{\bar R}  i^{2-\alpha}(\bar\rho) (v^2 + v_x^2) dx
\les \omega \mathscr{D}_1(t) + \omega^{-1} \mathcal{D}_1(t), \\
\mathcal{L}_{32}^{II} & \les \int_{\bar R/2}^{\bar R}  (x|v_x|+|v|)\int_0^x i^{-\alpha}(\bar\rho)y\lt(y|v_y|+|v|\rt)   dydx\\
&\les \int_{\bar R/2}^{\bar R}  (x|v_x|+|v|)\lt( \theta^{-1}   \mathscr{D}_1(t)  \rt)^{1/2}  dx\\
& \les \omega \mathscr{D}_1(t) + \omega^{-1} \theta^{-1}\mathcal{D}_1(t)
\end{align*}
for any $\omega>0$. Thus, one has that for any $\omega>0$,
\begin{align}
|\mathcal{L}_{3}| \les \omega \mathscr{D}_1(t) + (1+\omega^{-1})\theta^{-1} \mathcal{D}_1(t).\label{k3}
\end{align}

So, \eqref{tzeta0} follows form \eqref{zetakkk}-\eqref{k3}.

{\em Step 3}.
In this step, we prove that  for any $\omega>0$,
\begin{align}\label{tka-1}
&(1+t)^{2\zeta-1}\mathcal{E}_0(t) +(1+t)^{2\zeta}\left(\mathcal{E}_1
+\mathcal{D}_0\right)(t)
+\int_0^t \lt\{(1+s)^{2\zeta-1}\mathcal{D}_0(s)\rt.
\notag\\
&\lt.
+(1+s)^{2\zeta}\mathcal{D}_1(s)
\rt\}ds\les(\mathcal{E}_0+\mathcal{E}_1+\mathcal{D}_0)(0)
\notag
\\&+\omega \int_0^t  \mathscr{D}_0(s)ds
+\omega^{1-q_3}\int_0^t
(1+s)^{q_3(2\zeta-2)}
\mathscr{E}_0(s)ds,
\end{align}
where $q_3=2\bar\gamma/(2+2\bar\gamma\theta-\theta)$ satisfies $q_3(2\zeta-2)<-1$.
A suitable combination of  $\int_0^t(1+s)^{2\zeta-1}\eqref{r3x3}ds$  with  \eqref{l2tl} for $\ell=2\zeta$ gives that
\begin{align}\label{tka-2}
&(1+t)^{2\zeta-1}\mathcal{E}_0(t) +(1+t)^{2\zeta}\left(\mathcal{E}_1
+\mathcal{D}_0\right)(t)
\notag\\
&+\int_0^t \lt\{(1+s)^{2\zeta-1}\mathcal{D}_0(s)
+(1+s)^{2\zeta}\mathcal{D}_1(s)
\rt\}ds\notag\\
\les& (\mathcal{E}_0+\mathcal{E}_1
+\mathcal{D}_0)(0)
+(1+t)^{2\zeta-1}\mathcal{E}_1(t)
\notag\\
& +\int_0^t \lt\{(1+s)^{2\zeta-2}
\mathcal{E}_0(s)
+(1+s)^{2\zeta-1}\mathcal{D}_1(s)\rt\}ds\notag\\
\les & (\mathcal{E}_0+\mathcal{E}_1+\mathcal{D}_0)(0)+\int_0^t(1+s)^{2\zeta-2}
\mathcal{E}_0(s)ds,
\end{align}
where $2\zeta-1<1$ and \eqref{pr2v} have been used to derive the last inequality.
We choose   constants $q_3$, $q_4$ and $\delta_4$ satisfying
$$ \frac{1}{q_3}=\frac{1}{\bar\gamma}+\theta
-\frac{\theta}{2\bar\gamma}, \ \  \frac{1}{q_4}=1-\frac{1}{q_3}
 \ \ {\rm and}  \  \ \delta_4=\frac{1}{\alpha q_4-1},$$
then
$q_3,q_4>1$,       $q_3 (2\zeta-2) <-1$   and $0<\delta_4<\bar\gamma-1$.
In view of \eqref{irho} and \eqref{April4}, we see that
\begin{align*}
& i^{\alpha q_4/q_3}(\bar\rho)= i^{-\alpha }(\bar\rho) i^{\alpha q_4}(\bar\rho)\le 4 i^{-\alpha }(\bar\rho) p(\bar\rho) \bar\rho^{-1} i^{\alpha q_4-1}(\bar\rho) \\
&\les i^{-\alpha }(\bar\rho)  p(\bar\rho)  \bar\rho^{\delta_4 (\alpha q_4-1)-1} =i^{-\alpha }(\bar\rho) p(\bar\rho),
\end{align*}
which, together with the H\"{o}lder inequality, implies
\begin{align*}
\mathcal{E}_0\les \mathscr{E}_0^{1/q_3}
\lt(\int  i^{\alpha q_4/q_3}(\bar\rho) (|r-x|^2+|xr_x-x|^2)dx \rt)^{1/q_4}
\les \mathscr{E}_0^{1/q_3} \mathscr{D}_0^{1/q_4}.
\end{align*}
Substitute this into \eqref{tka-2} and use the Young inequality to prove \eqref{tka-1}.

{\em Step 4}. This step devotes to proving \eqref{12-21}.
It follows from \eqref{vvt}, and a suitable combination of \eqref{imt0}, \eqref{tzeta0}  and  \eqref{tka-1} with small $\omega$ that
\begin{align*}
&\mathscr{E}_0(t)+(1+t)^\zeta\mathscr{D}_0(t)
+(1+t)^{2\zeta}\left(\mathcal{E}_1
+\mathcal{D}_0\right)(t)
+(1+t)^{2\zeta-1}\mathcal{E}_0(t)\notag\\
&+\int_0^t\lt[\mathscr{D}_0(s)
+(1+s)^\zeta\mathscr{D}_1(s)
+(1+s)^{2\zeta}\mathcal{D}_1
(s)+
(1+s)^{2\zeta-1}\mathcal{D}_0(s)\rt]ds\notag\\
& \les  \lt(\mathscr{E}_0+ \mathscr{D}_0\rt) (0)
+
\theta^{-1}(\mathcal{E}_0+\mathcal{D}_1+\mathcal{E}_2)(0) \\
& +\int_0^t
 \lt\{  (1+s)^{- q_1 \zeta}+ \theta^{-1}(1+s)^{-2\zeta}+ (1+s)^{-q_3(2-2\zeta)} \rt\}
\mathscr{E}_0(s)ds,
\end{align*}
which, together with  $q_1 \zeta>1$, $2\zeta>1$, $q_3(2-2\zeta)>1$ and the Gronwall inequality, implies that
\begin{align}\label{zetadec-1}
&\mathscr{E}_0(t)+(1+t)^\zeta\mathscr{D}_0(t)
+(1+t)^{2\zeta}\left(\mathcal{E}_1
+\mathcal{D}_0\right)(t)
+(1+t)^{2\zeta-1}\mathcal{E}_0(t)\notag\\
&+\int_0^t\lt[\mathscr{D}_0(s)
+(1+s)^\zeta\mathscr{D}_1(s)
+(1+s)^{2\zeta}\mathcal{D}_1
(s)+
(1+s)^{2\zeta-1}\mathcal{D}_0(s)\rt]ds\notag\\
& \les  C(\theta^{-1}) \lt(\mathscr{E}_0 +\mathcal{D}_1+\mathcal{E}_2\rt)(0)
.
\end{align}
Based on \eqref{zetadec-1} and \eqref{vvt}, we integrate the product of $(1+t)^{2\zeta}$ and \eqref{phivt} over $[0,t]$ to obtain
\begin{align}
&(1+t)^{2\zeta}\mathcal{E}_2(t)+(1+t)^{2\zeta}\int  p(\bar\rho) (v^2 +x^2v_x^2)dx
+\int_0^t (1+s)^{2\zeta} \mathcal{D}_2(s)ds \notag\\
& \les C(\theta^{-1}) \lt(\mathscr{E}_0 +\mathcal{D}_1+\mathcal{E}_2\rt)(0).
\label{l2decay}
\end{align}
Indeed, the derivation of \eqref{l2decay} is the same as that of \eqref{21-10-13-1}. In view of \eqref{12.13.a}, we see that
\begin{align*}
(1+t)^{2\zeta}\mathcal{D}_1(t) \les \mathcal{D}_1(0) + \int_0^t (1+s)^{2\zeta} (\mathcal{D}_1+\mathcal{D}_2)(s)ds  ,
\end{align*}
which, together with \eqref{zetadec-1} and \eqref{l2decay}, proves \eqref{l2palphaest}.

It follows from   \eqref{hardy1}, \eqref{hardy2} and \eqref{har-1} that
\begin{align}\label{12.15-1}
 \int (r-x)^2 dx
\les \int x^2 i^2(\bar\rho)(|r-x|^2+|r_x-1|^2) dx.
\end{align}
When $\bar\gamma\le 2$, we set
$$\frac{1}{q_5}=\frac{2-\bar\gamma}{\bar\gamma}+\theta, \ \ \frac{1}{q_6}=1-\frac{1}{q_5} \ {\rm and} \ \delta_5=\frac{1}{2q_6-1},$$
then
$q_5,q_6>1$ and $0<\delta_5<\bar\gamma-1$.
So, it follows from \eqref{irho} and \eqref{April4} that
$$
 i^{2 q_6}(\bar\rho)\le 4  p(\bar\rho) \bar\rho^{-1} i^{2 q_6-1}(\bar\rho) \les    p(\bar\rho)  \bar\rho^{\delta_5 (2 q_6-1)-1} = p(\bar\rho),
$$
which gives, with the aid of the H\"{o}lder inequality and \eqref{12.15-1}, that
\begin{align}
 \int (r-x)^2 dx
\les \mathcal{E}_0^{1/q_5}\lt(\int  i^{2 q_6}(\bar\rho) (|r-x|^2+|xr_x-x|^2)dx \rt)^{1/q_6}
\les \mathcal{E}_0^{1/q_5}\mathcal{D}_0^{1/q_6}.\label{12.21-2}
\end{align}
When $\bar\gamma > 2$, it follows from \eqref{d-2} and \eqref{d-3} that
$$p''(0)=\lim_{t\to 0+} \frac{p'(s)}{s} =\bar\gamma\lim_{t\to 0+} \frac{p(s)}{s^2} =0 .$$
Then, we use  the L'Hopital rule, \eqref{d-0}, \eqref{d-00} and \eqref{21-11-2} to get
$$
\lim_{s\to 0+}\frac{i^2(s)}{p(s)}= 2\lim_{s\to 0+}\frac{ i(s)}{s}=2\lim_{s\to 0+}\frac{p'(s)}{s}=2p''(0)=0,
$$
so that
$i^2(\bar\rho)\les p(\bar\rho)$   and $ \int (r-x)^2 dx
\les \mathcal{D}_0 $,
where \eqref{12.15-1} has been used to derive the last inequality. This, together with \eqref{12.21-2} and \eqref{l2palphaest}, proves \eqref{rxt00}.
\hfill$\Box$

\begin{lemma}\label{lempoint}
Let $\Upsilon$ be defined by
\eqref{22.1.11-1}. Then it holds that for
$\theta\in \lt(0, 1-5/(4\bar\gamma)\rt]$ and $(x,t)\in I\times [0, T]$,
\begin{subequations}\label{22.1.5-2}\begin{align}
&
\lt(\Upsilon (1+t)^{2- {2}/{\bar\gamma}- {3\theta}/{2}}
+(1-\Upsilon)(1+t)^{2\zeta-{1}/{2}}\rt)x|r(x, t)-x|^2 \notag \\
&\quad +(1+t)^{2\zeta}x v^2(x,t) \les C(\theta^{-1}) \lt(\mathscr{E}_0 +\mathcal{D}_1+\mathcal{E}_2\rt)(0) , \label{12.21-da}\\
& x^3|r_x(x, t)-1|^2  \les C(\theta^{-1}) \lt(\mathscr{E}_0 +\mathcal{D}_1+\mathcal{E}_2\rt)(0)+x^3|r_{0x}(x) -1|^2,\label{12.21-1}
\end{align}\end{subequations}
where $\alpha=1-\theta$ and $\zeta=1- {1}/({2\bar\gamma})
- \theta/2$.
\end{lemma}

{\em Proof}. The bound for $xv^2(x, t)$ follows from \eqref{l2palphaest} and
\begin{align*}
& xv^2(x, t)=\int_0^x \lt(y v^2(y, t)\rt)_y dy \le \int v^2(y, t)dy \notag\\
&+2\lt(\int v^2(y, t) dy\rt)^{1/2}\lt(\int y^2 v^2_y(y,t) dy\rt)^{1/2} \le  3\mathcal{D}_1(t).
\end{align*}
Similarly, the bound for $x|r(x, t)-x|^2 $ follows from \eqref{12-21}. This finishes the proof of \eqref{12.21-da}.

To bound  $x^3|r_x(x, t)-1|^2 $,
we integrate   \eqref{remaineq} over $[x,\bar R]$  to obtain
\begin{align}\label{Zp}
\nu Z_t+ p\lt(\frac{x^2 \bar\rho}{r^2}
\rt)- p\lt(\frac{x^2 \bar\rho}{r^2 r_x}
\rt)
 = \mathscr{L}_1,
\end{align}
where
\begin{align*}
&Z(x,t)=\ln r_x,\\
&\mathscr{L}_1(x,t)=\lt( \frac{4}{3}\nu_1-2\nu_2\rt) \lt(\ln \frac{ r}{x}\rt)_t-\int_x^{\bar R} \mathscr{L}_2(y,t)dy,\\
&\mathscr{L}_2(x,t)=\bar\rho \frac{x^2}{r^2} v_t+\lt[p\lt(\frac{x^2 \bar\rho}{r^2}
\rt)\rt]_x-\frac{x^4}{r^4} p(\bar\rho)_x -4\nu_1 \lt(\frac{v}{r}\rt)_x.
\end{align*}
It follows from the Taylor expansion that there exist constants $\vartheta_4, \vartheta_5 \in (0,1) $ such that
\begin{align}
&Z(x,t)=\frac{1}{\vartheta_4 r_x+1-\vartheta_4} (r_x-1), \label{12.21-4}\\
&p\lt(\frac{x^2 \bar\rho}{r^2}
\rt)- p\lt(\frac{x^2 \bar\rho}{r^2 r_x}
\rt)=p'\lt(\frac{x^2 \bar\rho}{r^2 r_x}\lt(
\vartheta_5 r_x +1-\vartheta_5 \rt)\rt)\frac{x^2 \bar\rho}{r^2 r_x}(r_x-1). \notag
\end{align}
Thus,
\begin{align*}
&h(x,t)=\frac{1}{  p(\bar\rho) Z}\lt\{ p\lt(\frac{x^2 \bar\rho}{r^2}
\rt) -p\lt(\frac{x^2 \bar\rho}{r^2 r_x}
\rt)\rt\}\\
& =\frac{\vartheta_4 r_x+1-\vartheta_4}{\vartheta_5 r_x +1-\vartheta_5 } \frac{\check{\rho} p'(\check{\rho})}{  p(\bar\rho)}, \
{\rm where} \ \check{\rho}=\frac{x^2 \bar\rho}{r^2 r_x}\lt(
\vartheta_5 r_x +1-\vartheta_5 \rt).
\end{align*}
Due to \eqref{priori1} and the smallness of $\varepsilon_0$, one has
$$  \frac{1}{2} \frac{\check{\rho} p'(\check{\rho})}{  p(\bar\rho)} \le   h(x,t) \le 2 \frac{\check{\rho} p'(\check{\rho})}{ p(\bar\rho)} \ \  {\rm and} \  \ \frac{1}{2}  \bar\rho\le \check{\rho}\le \frac{3}{2} \bar\rho,$$
which, together with \eqref{21.3.15} and \eqref{til-2}, implies that
\begin{align}\label{12-17-a}
0<4^{-1} K_2 \le h(x,t) \le 4 K_3 < \infty.
\end{align}
So,   \eqref{Zp} can be rewritten  as
\begin{align}\label{12.21-h}
\nu Z_t+h(x,t)  p(\bar\rho) Z=\mathscr{L}_1,
\end{align}
which means
\begin{align*}
&Z(x, t)=\exp\lt\{-\nu^{-1}p(\bar\rho) \int_0^t h(x, \tau)d\tau \rt\}  Z(x,0)
\notag \\
& +\nu^{-1}\int_0^t \exp\lt\{-\nu^{-1}p(\bar\rho) \int_s^t h(x, \tau)d\tau \rt\}\mathscr{L}_1(x,s) ds.
\end{align*}
This gives,   with the aid of \eqref{12.21-4} and \eqref{12-17-a}, that
\begin{align}\label{zsim-1}
&\lt|r_x(x, t)-1\rt| \les \lt|r_{0x}(x)-1\rt|
+ \lt|\int_0^t \exp\lt\{-\nu^{-1}p(\bar\rho) \int_s^t h(x, \tau)d\tau \rt\}\mathscr{L}_1(x,s) ds\rt|.
\end{align}

For the first term in $\mathscr{L}_1$, we use the integration by parts, \eqref{12-17-a} and the Taylor expansion to get
\begin{align}
&\lt|\int_0^t \exp\lt\{-\nu^{-1}p(\bar\rho) \int_s^t h(x, \tau)d\tau \rt\} \lt(\ln \frac{ r(x,s)}{x}\rt)_s  ds\rt|\notag\\
& \les \sup_{s\in[0,t]}\ln \frac{ r(x,s)}{x} \les \frac{1}{x}  \sup_{s\in[0,t]} \lt| r(x,s)-x\rt|, \label{12.21-3}
\end{align}
by noticing that
\begin{align}
& p(\bar\rho)  \int_0^t \exp\lt\{-\nu^{-1}p(\bar\rho) \int_s^t h(x, \tau)d\tau \rt\}  ds\notag\\
& \le  p(\bar\rho)   \int_0^t  \exp\lt\{- \nu^{-1}p(\bar\rho) 4^{-1} K_2  (t-s)  \rt\}  ds
\le 4\nu K_2^{-1}. \label{12.20-a}
\end{align}
For the second term in $\mathscr{L}_1$, we have
\begin{align}
&\lt|\int_x^{\bar R} \mathscr{L}_2(y,t)dy\rt|
\les \int_x^{\bar R} y^{-2}\lt(y^2 \bar\rho|v_t|+ y|v_y|+|v|
\rt)dy\notag\\
&\quad+\lt|\frac{x^4}{r^4} p(\bar\rho) -p\lt(\frac{x^2 \bar\rho}{r^2}
\rt) \rt|+\int_x^{\bar R} y^{-2}p(\bar\rho)
\lt(y|r_y-1|+|r-y|\rt)dy\notag\\
&\les
x^{-\frac{3}{2}}\lt(\mathcal{E}_2^{\frac{1}{2}}+
\mathcal{D}_1^{\frac{1}{2}}\rt)(t)
+x^{-\frac{3}{2}}p(\bar\rho) \lt(x^{\frac{1}{2}}\lt|r(x,t)-x\rt|
+ \mathcal{E}_0^{\frac{1}{2}}(t)\rt),\label{12.19-a}
\end{align}
due to $ p(\bar\rho)\ge 0$, $(p(\bar\rho))_x\le 0$,  and  the following estimate:
\begin{align}
&\lt|\frac{x^4}{r^4} p(\bar\rho) -p\lt(\frac{x^2 \bar\rho}{r^2}
\rt) \rt|\le  p(\bar\rho) \lt|\frac{x^4}{r^4} -1 \rt|+\lt|p(\bar\rho) -p\lt(\frac{x^2 \bar\rho}{r^2}
\rt) \rt|\notag\\
&\quad\le p(\bar\rho) \lt|\frac{x^4}{r^4} -1 \rt|+ 4K_3 p(\bar\rho) \lt|1 -\frac{x^2 }{r^2}
\rt|.\label{12.19.1}
\end{align}
Here the derivation of \eqref{12.19.1} is the same as that of \eqref{12-17-a}.
In view of   \eqref{12.19-a}, \eqref{12-17-a}, \eqref{12.20-a}, $1<2\zeta$ and $\mathcal{E}_2\les \mathcal{D}_2$, we see that
\begin{align*}
& x^{{3}/{2}} \int_0^t \exp\lt\{-\nu^{-1}p(\bar\rho) \int_s^t h(x, \tau)d\tau \rt\}\lt|\int_x^{\bar R} \mathscr{L}_2(y,s)dy\rt| ds\notag\\
& \les \lt(\int_0^t (1+s)^{-2\zeta}ds\rt)^{1/2}\lt(\int_0^t (1+s)^{2\zeta}\lt(\mathcal{D}_1 +
\mathcal{D}_2 \rt)(s)ds\rt)^{1/2} \notag \\
 & +  p(\bar\rho) \int_0^t \exp\lt\{-\nu^{-1}p(\bar\rho) \int_s^t h(x, \tau)d\tau \rt\} ds \sup_{s\in [0,t]} \lt( x^{1/2}|r(x,s)-x|+ \mathcal{E}_0^{{1}/{2}}(s) \rt)\notag\\
& \les \lt(\int_0^t (1+s)^{2\zeta}\lt(\mathcal{D}_1 +
\mathcal{D}_2 \rt)(s) ds\rt)^{1/2}+  \sup_{s\in [0,t]} \lt( x^{1/2}|r(x,s)-x|+ \mathcal{E}_0^{{1}/{2}}(s) \rt).
\end{align*}
This, together with \eqref{zsim-1}, \eqref{12.21-3}, \eqref{12.21-da} and \eqref{l2palphaest}, proves \eqref{12.21-1}.
\hfill $\Box$

\subsection{Higher-order estimates}\label{22.1.4-3}
This subsection consists of Lemma \ref{22.1.4-4}.

\begin{lemma}\label{22.1.4-4}
Let $\Upsilon$ be defined by
\eqref{22.1.11-1} and $L^\infty=L^\infty(I)$. Then it holds that for $l\in (0, \bar R)$,
$\theta\in \lt(0, 1-5/(4\bar\gamma)\rt]$ and $t\in [0, T]$,
\begin{subequations}\label{22.1.5-1}
\begin{align}
&\mathfrak{E}(t)+ \lt\| (D^{0,0}+D^{1,0})(\cdot,t)\rt\|_{L^\infty}\les \mathfrak{E}(0), \label{22.1.5}\\
&(1+t)^{2\zeta-1/2} \lt\{  (1-\Upsilon) \lt\| (r-x)^2(\cdot,t)\rt\|_{L^\infty}  +\lt\|v^2(\cdot,t)\rt\|_{L^\infty}  \rt\} \notag\\
  & \ \ +(1+t)^{\zeta} \lt\| (x^3v_x^2)(\cdot,t)\rt\|_{L^\infty}  +
 \Upsilon(1+t)^{2- {2}/{\bar\gamma}-3\theta/2} \lt\| (r-x)^2(\cdot,t)\rt\|_{L^\infty}
\notag\\
& \ \ +(1+t)^{2\zeta-1} \lt\{ \lt\| {D}^{1,0} (\cdot,t)\rt\|_{L^\infty}
+\Upsilon\lt\| \lt(\bar\rho^{(3\bar\gamma-2)/2} \mathcal{Q}^2\rt)(\cdot,t)\rt\|_{L^\infty} \rt.\notag\\
& \ \ \lt. +\int \lt(\bar\rho v^2_t + {D}^{0,0} \rt)  (x, t) dx  \rt\} \les C(\theta^{-1})\mathfrak{E}(0),\label{12.28-9}\\
&(1+t)^{\min\lt\{\frac{4-\bar\gamma}{2\bar\gamma}
-\theta\lt(\frac{\bar\gamma-1}{2\bar\gamma}+ \frac{4-\bar\gamma-\theta(\bar\gamma-1)}
{4\bar\gamma-2}\rt), \ 2\zeta-1 \rt\}} \lt\| (\bar\rho^{2}\mathcal{Q}^2)
(\cdot,t)\rt\|_{L^\infty}\notag\\
& \ \ \les C(\theta^{-1})\mathfrak{E}(0),
\ \ {\rm if} \ \ 2 < \bar\gamma<4 \ \ {\rm and}\  \ \theta<({4-\bar\gamma})/({\bar\gamma-1}),\label{20220111} \\
&(1+t)^{2\zeta-1} \lt\{ \int_0^{\bar R -l} \lt({D}^{0,1}+{D}^{1,1}\rt)(x,t) dx \rt.
\notag\\
& \ \ \lt. +\lt\| {D}^{0,0} (\cdot, t)\rt\|_{L^\infty\lt([0, \bar R-l]\rt)} \rt\}
 \les C\lt(\theta^{-1}, \ l^{-1}\rt)\mathfrak{E}(0),\label{22.1.6-1}
\end{align}
\end{subequations}
where $\zeta=1- {1}/({2\bar\gamma})
- \theta/2$,
$
{D}^{i,j}  = |\partial_t^i\partial_x^j (r/x-1)|^2
+|\partial_t^i\partial_x^j(r_{x}-1)|^2
$ and
$\mathcal{Q}={x^2}/({r^2r_x})-1$.
\end{lemma}

{\em Proof}. The proof of this lemma is based on the following estimates:
\begin{subequations}\label{12.26-1}\begin{align}
&\int w\lt(10|(r/x)_x|^2+|r_{xx}|^2\rt)dx\le 2 \int w |\mathcal{Q}_x|^2 dx, \label{wrxx}\\
&\int w \lt(10|(v/x)_x|^2+|v_{xx}|^2\rt)dx\le  6\int w\lt(|\mathcal{Q}_{xt}|^2+32 |\mathcal{Q}_x|^2 \rt)dx,\label{wvxx}
\end{align}\end{subequations}
for any given  function $w(x)$ on $[0, \bar R]$ satisfying $w(x)\geq 0$ and $w'(x)\leq 0$ on $[0, \bar R],$  and $w(\bar R)=0$, where
\begin{equation}\label{12.25-1}
\mathcal{Q}={x^2}/({r^2r_x})-1.
\end{equation}
Indeed, the proof of \eqref{12.26-1} is the same as that of Lemma 3.9 in \cite{LZeng}, so we omit the detail here.

{\em Step 1}. In this step, we prove that
\begin{subequations}\begin{align}
&\int \bar\rho^{-1} p^2(\bar\rho)  {D}^{0,1} dx + \int_0^t \int \bar\rho^{-1} p^3(\bar\rho)  {D}^{0,1}dx ds  \les \mathfrak{E}(0), \label{higlobal}\\
&(1+t)^{\zeta}x^3v_x^2(x, t)  \les  C(\theta^{-1})\mathfrak{E}(0), \label{22.1.6-2} \\
& \Upsilon(1+t)^{2\zeta-1} x^3\bar\rho^{(3\bar\gamma-2)/2}(x)\mathcal{Q}^2(x,t)
\les C(\theta^{-1})\mathfrak{E}(0),  \label{22.1.10-2}
\end{align}\end{subequations}
and for $ 2 < \bar\gamma<4$ and $\theta<({4-\bar\gamma})/({\bar\gamma-1})$,
\begin{align} (1+t)^{\min\lt\{\frac{4-\bar\gamma}{2\bar\gamma}
-\theta\lt(\frac{\bar\gamma-1}{2\bar\gamma}+ \frac{4-\bar\gamma-\theta(\bar\gamma-1)}
{4\bar\gamma-2}\rt), \ 2\zeta-1 \rt\}} x^3\bar\rho^{2}(x)\mathcal{Q}^2(x,t)
\les C(\theta^{-1})\mathfrak{E}(0). \label{22.1.10-4}
\end{align}

We  use \eqref{pxphi} to
rewrite \eqref{fixedp-a}  as
\begin{align}\label{remainq}
&\nu \lt((1+ \mathcal{Q})^{-1}  \mathcal{Q}_x\rt)_t+  \bar\rho  p'\lt( \bar\rho(1+\mathcal{Q})\rt)\mathcal{Q}_{x}=\mathcal{H},
\end{align}
where $\mathcal{Q}$ is defined by \eqref{12.25-1}, and
\begin{align*}
& \mathcal{H} =
\lt(\frac{x^4}{r^4}-1-\mathcal{Q}\rt) \lt(p(\bar\rho)\rt)_x
+\lt(p'(\bar\rho)- p'\lt(\bar\rho(1+\mathcal{Q})\rt)\rt) (1+\mathcal{Q}) \bar\rho_x-\bar\rho \frac{x^2}{r^2} v_t.
\end{align*}
In view of the Taylor expansion, \eqref{pxphi}, \eqref{priori1},  \eqref{21.3.15} and \eqref{til-2},
we see that
\begin{subequations}\label{12.25-2}
\begin{align}
&4^{-1}K_2 p(\bar\rho)\le \bar\rho  p'\lt( \bar\rho(1+\mathcal{Q})\rt)\le 4K_3p(\bar\rho), \label{12.25-2.a} \\
&|\mathcal{H}|\les  \bar\rho \lt(x|r_x-1|+|r-x|+|v_t|\rt).
\end{align}
\end{subequations}
In fact, the second term in $\mathcal{H}$ can be bounded as follows:
\begin{align*}
&\lt|\lt( p'\lt(\bar\rho(1+\mathcal{Q})\rt)-p'(\bar\rho) \rt)   \bar\rho_x \rt|=
\frac{\bar\rho \lt|p''\lt(\bar\rho \lt(1+\vartheta_6\mathcal{Q}\rt)\rt) \rt|}{p'(\bar\rho)} x\bar\rho\phi\lt|\mathcal{Q}\rt|
\\
& \les \frac{  p'\lt(\bar\rho \lt(1+\vartheta_6\mathcal{Q}\rt)\rt) }{p'(\bar\rho)} {x\bar\rho } \lt|\mathcal{Q}\rt|
\les \frac{  p\lt(\bar\rho \lt(1+\vartheta_6\mathcal{Q}\rt)\rt) }{p(\bar\rho)} {x\bar\rho }\lt|\mathcal{Q}\rt|
 \les x\bar\rho \lt|\mathcal{Q}\rt|
\end{align*}
for some constant  $\vartheta_6\in (0,1)$.
We integrate the product of  \eqref{remainq} and $\bar\rho^{-1} p^2(\bar\rho) (1+ \mathcal{Q})^{-1}  \mathcal{Q}_x$ over $I$, and use the Cauchy inequality and \eqref{12.25-2} to obtain
\begin{align}\label{qx}
&\frac{\nu}{2}\frac{d}{dt}\int \bar\rho^{-1} p^2(\bar\rho)(1+\mathcal{Q})^{-2} \mathcal{Q}_{x}^2 dx +
\frac{K_2}{8} \int \bar\rho^{-1} p^3(\bar\rho)  \mathcal{Q}_{x}^2 dx
\notag \\
& \les \int \bar\rho p(\bar\rho) \lt(v_t^2 + x^2|r_x-1|^2+|r-x|^2\rt),
\end{align}
which gives, with the aid of   \eqref{vvt} and the fact that $f(x,t)-f(0,t)=\int_0^x f_y(y,t)dy$ for any function $f$, that
\begin{align}\label{wl2qx}
\int \bar\rho^{-1} p^2(\bar\rho)  \mathcal{Q}_{x}^2 dx +
\int_0^t \int \bar\rho^{-1} p^3(\bar\rho)  \mathcal{Q}_{x}^2 dxds
\les \mathfrak{E}(0).
\end{align}
Due to \eqref{pxphi}, \eqref{pcon2} and \eqref{d-3}, one has  that for any constant $k > 3/4$,
\begin{subequations}\label{12.25-3}\begin{align}
&  \bar\rho^{-1} p^k(\bar\rho)
=0 \ {\rm at} \ x=\bar R,\\
&\lt( \bar\rho^{-1} p^k(\bar\rho) \rt)_x
=-\frac{x \phi}{\bar\rho p'(\bar\rho)}p^{k-1}(\bar\rho)
\lt(k\bar\rho p'(\bar\rho)-p(\bar\rho)\rt)\le 0,
\end{align}\end{subequations}
which, together with \eqref{wl2qx} and \eqref{wrxx}, proves \eqref{higlobal}.

To prove \eqref{22.1.6-2}, we notice that
\begin{align}
&\lt|v_x(x,t)\rt|=\lt|r_xZ_t\rt|=\nu^{-1} r_x \lt|\mathscr{L}_1- h p(\bar\rho) Z\rt|\les \lt|\mathscr{L}_1\rt|+  p(\bar\rho)\lt|  r_x-1\rt|\notag\\
& \les \lt|\int_x^{\bar R} \mathscr{L}_2(y,t)dy\rt| +x^{-1}|v(x,t)|  + p(\bar\rho)\lt|  r_x(x,t)-1\rt|,\notag
 \end{align}
which is due to \eqref{12.21-4}-\eqref{12.21-h};
\begin{align*}
\lt|\int_x^{\bar R} \mathscr{L}_2(y,t)dy\rt|
 \les
x^{-\frac{3}{2}}\lt(\mathcal{E}_2^{\frac{1}{2}}+
\mathcal{D}_1^{\frac{1}{2}}+\mathcal{D}_0^{\frac{1}{2}} \rt)(t)
+x^{-1} p(\bar\rho) \lt|r(x,t)-x\rt| ,
\end{align*}
which follows from the same derivation as that of \eqref{12.19-a}; and
 \begin{align*}
&x^3 p^2(\bar\rho(x))\lt(r_x(x,t)-1\rt)^2 =\int_0^x
\lt( y^3 p^2(\bar\rho(y))\lt(r_y(y,t)-1\rt)^2\rt)_y dy\\
& \ \les \mathcal{D}_0(t)+  \int y^3 p^2(\bar\rho(y))|r_y(y,t)-1||r_{yy}(y,t)| dy \\
& \  \les \mathcal{D}_0(t)+
\lt( \int \bar\rho^{-1}
 p^2(\bar\rho(y))  r_{yy}^2(y,t)  dy \rt)^{1/2}  \mathcal{D}_0^{1/2}(t),\\
& x  p^2(\bar\rho(x))\lt(r(x,t)-x\rt)^2 \les \mathcal{D}_0(t),
\end{align*}
which is due to \eqref{pxphi}.
Then, \eqref{22.1.6-2} follows from \eqref{higlobal}, \eqref{12.21-da} and \eqref{l2palphaest}.

When $\bar\gamma\le 2$, we use
 \eqref{pxphi} and the  H$\ddot{o}$lder inequality to get
\begin{align}
&x^3\bar\rho^{(3\bar\gamma-2)/2}\mathcal{Q}^2= \int_0^x \lt(y^3 \bar\rho^{(3\bar\gamma-2)/2} \mathcal{Q}^2 \rt)_y dy
\notag\\
& \le  3 \int_0^x  y^2 \bar\rho^{(3\bar\gamma-2)/2} \mathcal{Q}^2 dy+2\int_0^x  y^3 \bar\rho^{(3\bar\gamma-2)/2} \mathcal{Q} \mathcal{Q}_y dy
\notag\\
& \les  \mathcal{E}_0   + \lt(\int \bar\rho^{-1} p^2(\bar\rho)  \mathcal{Q}_{y}^2 dy \rt)^{1/2}
\lt(\int y^2 \bar\rho^{3\bar\gamma-1} p^{-2}(\bar\rho)  \mathcal{Q}^2 dy \rt)^{1/2}. \label{22.1.10-1}
\end{align}
We set
$$\frac{1}{q_7}=2\zeta-1, \ \ \frac{1}{q_8}=1-\frac{1}{q_7}
\ \ {\rm and} \ \
\delta_6=\frac{3q_8-1}{1+(3\bar\gamma-4)q_8},
$$
then $q_7,q_8>1$ and $1/\delta_6 >{\bar\gamma-1}$, which means, with the help of \eqref{irho} and \eqref{April4}, that
\begin{align*}
\bar\rho^{(3\bar\gamma-1)q_8}
p^{1-3q_8}(\bar\rho)
\les  \bar\rho^{1+(3\bar\gamma-4)q_8}
i^{1-3q_8}(\bar\rho)
\les i^{\delta_6(1+(3\bar\gamma-4)q_8)+1-3q_8}
(\bar\rho)=1.
\end{align*}
So, it follows from the  H$\ddot{o}$lder inequality that
\begin{align*}
&\int y^2 \bar\rho^{3\bar\gamma-1} p^{-2}(\bar\rho)  \mathcal{Q}^2 dy\le \lt(\int y^2 p(\bar\rho)    \mathcal{Q}^2 dy  \rt)^{1/q_7} \notag \\
&
\times\lt(\int y^2 \bar\rho^{(3\bar\gamma-1)q_8}p^{1-3q_8}(\bar\rho)    \mathcal{Q}^2 dy  \rt)^{1/q_8}
\le \mathcal{D}_0^{1/q_7} \mathcal{E}_0^{1/q_8},
\end{align*}
which, together with \eqref{22.1.10-1}, \eqref{wl2qx} and \eqref{l2palphaest}, proves \eqref{22.1.10-2}.

In a similar way to deriving \eqref{22.1.10-1}, we have
\begin{align}
x^3\bar\rho^2\mathcal{Q}^2
 \les  \mathcal{E}_0   + \lt(\int \bar\rho^{-1} p^2(\bar\rho)  \mathcal{Q}_{y}^2 dy \rt)^{1/2}
\lt(\int y^2 \bar\rho^{5} p^{-2}(\bar\rho)  \mathcal{Q}^2 dy \rt)^{1/2}. \label{22.1.10-3}
\end{align}
When $ 2<\bar\gamma<4$ and $\theta<({4-\bar\gamma})/({\bar\gamma-1})$, we set
$$\frac{1}{q_9}=\frac{4-\bar\gamma-\theta(\bar\gamma-1)}
{2\bar\gamma-1},  \ \frac{1}{q_{10}}=1-\frac{1}{q_9} \ {\rm and} \  \delta_7=\frac{(2-\alpha) q_9+1+\alpha}{3q_9-1},
$$
then $q_9,q_{10}>1$ and $1/\delta_7 >{\bar\gamma-1}$. So,
$$\bar\rho^{5q_9}p^{-2 q_9-1}(\bar\rho) i^{\alpha(q_9-1)}
\les  \bar\rho^{3q_9-1}i^{(\alpha-2) q_9-1-\alpha}(\bar\rho)
\les i^{\delta_7(3q_9-1)+(\alpha-2) q_9-1-\alpha }(\bar\rho)=1,$$
which means, using the  H$\ddot{o}$lder inequality, that
\begin{align}
&\int y^2 \bar\rho^{5} p^{-2}(\bar\rho)  \mathcal{Q}^2 dy \le
\lt(\int \bar\rho^{5q_9}p^{-2 q_9-1}(\bar\rho) i^{\alpha(q_9-1)} (\bar\rho)  y^2 p(\bar\rho) \mathcal{Q}^2 dy\rt)^{1/q_{9}}
\notag\\
&\times\lt(\int y^2 i^{-\alpha}(\bar\rho) \mathcal{Q}^2 dy\rt)^{1/q_{10}} \le \mathcal{D}_0^{1/q_9}\mathscr{E}_0^{1/q_{10}}.\notag
\end{align}
This proves \eqref{22.1.10-4}, by noting \eqref{22.1.10-3}, \eqref{wl2qx} and \eqref{l2palphaest}.

{\em Step 2}. In this step, we prove that
\begin{align}
(1+s)^{2\zeta-1}\int \bar\rho v^2_t  dx+\int_0^t (1+s)^{2\zeta-1}\int {D}^{2,0} dxds\les C(\theta^{-1})\mathfrak{E}(0)\label{glovt}.
\end{align}

In a similar way to deriving \eqref{qx},
one has
\begin{align}
&\frac{\nu}{2}\frac{d}{dt}\int \psi \bar\rho^{-1} p^2(\bar\rho)(1+\mathcal{Q})^{-2} \mathcal{Q}_{x}^2 dx +
\frac{K_2}{8} \int \psi \bar\rho^{-1} p^3(\bar\rho)  \mathcal{Q}_{x}^2 dx
\notag \\
& \les \int \psi \bar\rho p(\bar\rho) \lt(v_t^2 + x^2|r_x-1|^2+|r-x|^2\rt),
\label{psil2}
\end{align}
where $\psi$ is a
smooth cut-off function on $[0, \bar R]$ satisfying that for any fixed constant $l\in (0, \bar R)$,
\begin{equation}\label{psir}\begin{split}
&\psi=1 \ \ {\rm on} \  \ [0, \bar R-l], \ \ \psi=0 \ \ {\rm on} \ \  [\bar R-l/2, \bar R],\\
&{\rm and} \ \  -8/l\leq \psi'(x)\leq 0 \ \ {\rm and} \  \  [0, \bar R].
\end{split}
\end{equation}
In view of \eqref{irho}, \eqref{April4} and \eqref{har-1}, we see that on $[0,\bar R-l/2]$,
\begin{align}\label{12.27-1}
p^{-1}(\bar\rho)\les
\bar\rho^{-1} i^{-1}(\bar\rho)\les i^{-1-8/(7\bar\gamma-7)} \les (\bar R -x)^{-1-8/(7\bar\gamma-7)} \les C(l^{-1}),
\end{align}
which, together with \eqref{psil2}, \eqref{l2palphaest}, \eqref{wl2qx} and $2\zeta-1<1$, implies that
\begin{align}
&(1+t)^{2\zeta-1}\int \psi \bar\rho^{-1} p^2(\bar\rho) \mathcal{Q}_{x}^2 dx +\int_0^t
(1+s)^{2\zeta-1}
\int \psi \bar\rho^{-1} p^3(\bar\rho)  \mathcal{Q}_{x}^2 dxds
\notag \\
& \les C(\theta^{-1})\mathfrak{E}(0)
+\int_0^t
\int_0^{\bar R -l/2} \psi \bar\rho^{-1} p^2(\bar\rho)  \mathcal{Q}_{x}^2 dxds\notag\\
& \les C(\theta^{-1})\mathfrak{E}(0)
+C(l^{-1})\int_0^t
\int   \bar\rho^{-1} p^3(\bar\rho)  \mathcal{Q}_{x}^2 dxds\les C(\theta^{-1}, l^{-1})\mathfrak{E}(0).
\label{qxdec}
\end{align}
It follows from \eqref{remainq}, \eqref{21apri} and \eqref{12.25-2} that
\begin{align}
& |{Q}_{xt}| \les |\mathcal{Q}_{x}|+\bar\rho \lt(x|r_x-1|+|r-x|+|v_t|\rt),
\label{12.27-h}
\\
&\int \psi \bar\rho^{-1} p^3(\bar\rho)  {Q}_{xt}^2 dx \les  \int \psi \bar\rho^{-1} p^3(\bar\rho)  \mathcal{Q}_{x}^2 dx
+\mathcal{D}_0+\mathcal{D}_2,\notag
\end{align}
which means,  with the help of \eqref{l2palphaest} and \eqref{qxdec}, that
$$
\int_0^t   (1+s)^{2\zeta-1}\int \psi \bar\rho^{-1} p^3(\bar\rho)  {Q}_{xs}^2 dx ds \les C(\theta^{-1}, l^{-1})\mathfrak{E}(0).
$$
This, together with \eqref{12.26-1},   \eqref{12.25-3}, \eqref{psir} and \eqref{qxdec},   gives that
\begin{align}
&(1+t)^{2\zeta-1}\int \psi \bar\rho^{-1} p^2(\bar\rho) {D}^{0,1}(x,t) dx \notag\\
& +\int_0^t
(1+s)^{2\zeta-1}
\int \psi \bar\rho^{-1} p^3(\bar\rho)  {D}^{1,1}(x,s)  dxds
\les C(\theta^{-1}, l^{-1})\mathfrak{E}(0).
\label{wrvxx}
\end{align}

Take $\partial_t$ on \eqref{remaineq} and use \eqref{pxphi} to yield
\begin{align*}
&\bar\rho \frac{x^2}{r^2}v_{tt}
-\nu\lt(2\frac{v_t}{r}+\frac{v_{tx}}{r_x}\rt)_x
=-\bar\rho \lt(\frac{x^2}{r^2}\rt)_tv_{t}
\\
& -\lt[p\lt(\frac{x^2 \bar\rho}{r^2 r_x}
\rt)\rt]_{tx}+\lt(\frac{x^4}{r^4} \rt)_t \lt(p(\bar\rho)\rt)_x
-\nu\lt(2\frac{v^2}{r^2}
+\frac{v_x^2}{r_x^2}\rt)_x.
\end{align*}
Let $\psi$ be defined by \eqref{psir} with $l=\bar R/2$, we
integrate the product of the equation above and $\psi v_t$ over $[0, \bar R]$, and use the boundary condition $ v_t(0, t)=\psi(\bar R)=0$, \eqref{pxphi},  \eqref{21apri}, \eqref{12.25-2.a} and the Cauchy inequality to obtain
\begin{align}
&\frac{1}{2}\frac{d}{dt}\int \psi \bar\rho  \frac{x^2}{r^2} v_t^2 dx +\frac{\nu}{2} \int  \psi \lt(r_x\frac{v_t^2}{r^2}
+\frac{v_{tx}^2}{r_x}\rt) dx \notag\\
&\les (\mathcal{D}_1+\mathcal{D}_2)(t) +\int {D}^{1,0}(x,t)dx  ,\label{intvt}
\end{align}
where we have used the following estimate:
$$
 \int  \lt(2\frac{v_t}{r}+\frac{v_{tx}}{r_x}\rt) (\psi v_t)_x dx=  \int  \psi \lt(r_x\frac{v_t^2}{r^2}
+\frac{v_{tx}^2}{r_x}\rt) dx
+ \int \psi'\lt(\frac{v_t^2}{r}
+\frac{v_tv_{tx}}{r_x}\rt) dx.
$$
Due to \eqref{hardy1}, one has
\begin{align}\label{12.28.1}
\int {D}^{1,0}dx \les
\sum_{j=0,1}\int_0^{\bar R/2} x^2 {D}^{1,j}dx + \int_{\bar R/2}^{\bar R} x^2 {D}^{1,0} dx \les \int_0^{\bar R/2}  {D}^{1,1}dx+ \mathcal{D}_1,
\end{align}
which means, using \eqref{l2palphaest}, and \eqref{wrvxx} with $l=\bar R/2$, that
\begin{align}
 \int_0^t (1+s)^{2\zeta-1}\int {D}^{1,0}(x,s)  dxds\les C(\theta^{-1}) \mathfrak{E}(0). \label{iivx}
\end{align}
It follows from \eqref{intvt}, \eqref{l2palphaest}, \eqref{iivx} and $0<2\zeta-1<1$ that
$$
(1+s)^{2\zeta-1}\int \psi \bar\rho v^2_t  dx+\int_0^t (1+s)^{2\zeta-1}\int \psi {D}^{2,0} dxds\les C(\theta^{-1})\mathfrak{E}(0),
$$
which proves \eqref{glovt}, by use of \eqref{l2palphaest}.

{\em Step 3}. This step devotes to  proving \eqref{22.1.5-1}. In view of \eqref{12.27-h}, \eqref{qxdec}, \eqref{glovt} and \eqref{l2palphaest}, we see that
$$
(1+t)^{2\zeta-1}\int \psi \bar\rho^{-1} p^2(\bar\rho) {Q}_{xt}^2 dx
\les C(\theta^{-1}, l^{-1})\mathfrak{E}(0),
$$
where $\psi$ is defined by \eqref{psir} for any fixed constant $l\in (0, \bar R)$. Thus, we use \eqref{wvxx},   \eqref{12.25-3}, \eqref{psir} and \eqref{qxdec} to get
$$
(1+t)^{2\zeta-1}\int \psi \bar\rho^{-1} p^2(\bar\rho) {D}^{1,1}(x,t) dx \les C(\theta^{-1}, l^{-1})\mathfrak{E}(0),
$$
which, together with
\eqref{wrvxx} and \eqref{12.27-1}, implies
\begin{align}
(1+t)^{2\zeta-1} \int_0^{\bar R -l} \lt({D}^{0,1}+{D}^{1,1}\rt)(x,t) dx  \les C\lt(\theta^{-1}, \ l^{-1}\rt)\mathfrak{E}(0).  \label{locdec}
\end{align}
In a similar way to deriving \eqref{12.28.1}, we have
$\int {D}^{0,0}dx \les \int_0^{\bar R/2}  {D}^{0,1}dx+ \mathcal{E}_0$,
which means, with the aid of \eqref{12.28.1}, \eqref{l2palphaest}, and \eqref{locdec} with $l=\bar R/2$, that
\begin{align}\label{22.1.6-4}
(1+t)^{2\zeta-1} \int \lt({D}^{0,0}+D^{1,0}\rt)(x,t)dx \les  C(\theta^{-1})\mathfrak{E}(0).
\end{align}
Due to \eqref{locdec}, \eqref{22.1.6-4} and $\|\cdot\|_{L^\infty}\les \|\cdot\|_{H^1} $, one has
\begin{align}\label{22.1.7-1}
(1+t)^{2\zeta-1} \lt\|\lt({D}^{0,0} + {D}^{1,0} \rt)(\cdot, t)\rt\|_{L^\infty\lt([0, \bar R-l]\rt)}
 \les C\lt(\theta^{-1}, \ l^{-1}\rt)\mathfrak{E}(0),
\end{align}
which proves \eqref{22.1.6-1}, using \eqref{locdec}. \eqref{20220111} follows from \eqref{22.1.10-4} and \eqref{22.1.7-1} with $l=\bar R/2$.

It follows from the fact that $\|f\|_{L^\infty}^2\le 2\|f\|_{L^2}\|f_x\|_{L^2}$ for any function $f$ satisfying $f(x=0)=0$, \eqref{12-21} and \eqref{22.1.6-4}
 that
\begin{align}
&  \lt( \Upsilon(1+t)^{2- {2}/{\bar\gamma}-3\theta/2} + (1-\Upsilon) (1+t)^{2\zeta-1/2} \rt)  |r(x,t)-x|^2 \notag\\
& +(1+t)^{2\zeta-1/2}v^2(x,t)\les C(\theta^{-1})\mathfrak{E}(0), \label{22-1-7-1}
\end{align}
which, together with
\eqref{22.1.7-1} with $l=\bar R/2$ and \eqref{22.1.6-2}, means
\begin{align}
(1+t)^{2\zeta-1} D^{1,0} (x, t)  \les  C(\theta^{-1})\mathfrak{E}(0). \label{22-1-7-2}
\end{align}
\eqref{12.28-9} follows from \eqref{22.1.6-2},  \eqref{22.1.10-2}, \eqref{glovt}, and \eqref{22.1.6-4}-\eqref{22-1-7-2} with $l=\bar R/2$.
Letting $\theta=1/32$ and $l=\bar R/2$ in \eqref{12.21-1} and \eqref{22.1.7-1}-\eqref{22-1-7-2}  gives
$ \lt( D^{0,0}+{D}^{1,0} \rt)(x,t) \les \mathfrak{E}(0)$,
which proves \eqref{22.1.5}, by use of \eqref{glovt} with $\theta=1/32$ and \eqref{higlobal}.
\hfill $\Box$

\subsection*{Acknowledgements.}
This research was supported in part by National Sciences Fundation of China (NSFC)   Grants 12171267, 11822107 and 12101350;  China Postdoctoral Science Foundation
under grant 2021M691818; and a grant from the Research Grants Council of the Hong Kong Special Administrative Region, China
(Project No. 11307420).

\bibliographystyle{plain}

\end{document}